\documentclass[12pt]{article}
\usepackage[left=2cm,top=2cm,right=2cm,bottom=2cm]{geometry}

\usepackage[cp1251]{inputenc}
\usepackage[british]{babel}

\usepackage{amsfonts,amsmath,amssymb,graphicx,latexsym,xcolor,mathtools}
\usepackage{hyperref,enumitem,enumerate,float}
\usepackage{amsthm,verbatim}
\usepackage{cite}
\usepackage[colorinlistoftodos]{todonotes}

\numberwithin{equation}{section}
\pdfoutput=1

\theoremstyle{plain}
\newtheorem{theorem}{Theorem}[section]

\newtheorem{proposition}[theorem]{Proposition}

\theoremstyle{definition}

\newenvironment{rema}
{%
	\pushQED{\qed}\begin{rema/}}
	{\popQED\end{rema/}}
\newtheorem{defi/}[theorem]{Definition}
\newtheorem{rema/}[theorem]{Remark}
\newtheorem{exa/}[theorem]{Example}

\newcommand{\sgn}{\text{sgn}}

\newcommand{\ptl}{\partial}

\newcommand{\bdG}{\boldsymbol{\Gamma}}

\newcommand{\bdx}{\boldsymbol{x}}

\newcommand{\bdxi}{\boldsymbol{\xi}}
\newcommand{\bdomega}{\boldsymbol{\omega}}

\newcommand{\STOP}{{\color{red} STOPPED HERE}}
\makeatletter
\newcommand*{\transpose}{%
	{\mathpalette\@transpose{}}%
}
\newcommand*{\@transpose}[2]{%
	\raisebox{\depth}{$\m@th#1\intercal$}%
}
\makeatother
\newcommand{\vph}{\varphi}

\newcommand{\infixand}{\text{ and }}
\newcommand{\mathd}{\mathrm{d}}
\newcommand{\nobracket}{}
\newcommand{\tmcolor}[2]{{\color{#1}{#2}}}
\newcommand{\tmem}[1]{{\em #1\/}}
\newcommand{\tmmathbf}[1]{\ensuremath{\boldsymbol{#1}}}
\newcommand{\tmop}[1]{\ensuremath{\operatorname{#1}}}

\definecolor{myred}{RGB}{160,0,0}
\definecolor{mygreen}{RGB}{0,160,0}
\definecolor{myblue}{RGB}{0,0,160}
\hypersetup{
	colorlinks = true,
	linkcolor = {myred},
	citecolor = {mygreen},
	urlcolor  = {myblue},
}
\title{A contribution to the mathematical
	theory of diffraction\\
	Part II. Recovering the far-field asymptotics of the quarter-plane problem.}

\author{Rapha\"{e}l C. Assier$^{*}$, Andrey V. Shanin$^{\dagger}$ and Andrey I. Korolkov$^{*}$\\
	\footnotesize{$^{*}$ Department of Mathematics, University of Manchester, Oxford Road, Manchester, {\rm M13 9PL}, UK}\\
	\footnotesize{$^{\dagger}$ Department of Physics (Acoustics Division), Moscow State University, Leninskie Gory, {\rm 119992}, Moscow, Russia}
}

\begin{document}

\maketitle

\begin{abstract}
  We apply the stationary phase method developed in (Assier, Shanin \& Korolkov, QJMAM, 76(1), 2022) to the problem of wave diffraction by a quarter-plane. The wave field is written as a double Fourier transform of an unknown spectral function. We make use of the analytical continuation results of (Assier \& Shanin, QJMAM, 72(1), 2018) to uncover the singularity structure of this spectral function.  This allows us to provide a closed-form far-field asymptotic expansion of the field by estimating the  double Fourier integral near some special points of the spectral function. All the known results on the far-field asymptotics of the quarter-plane problem are recovered, and new mathematical expressions are derived for the secondary diffracted waves in the plane of the scatterer.    
\end{abstract}

\section{Introduction}

 The present article (Part II) is  a continuation of \cite{Part6A} (Part I), in which a general mathematical framework was developed for the asymptotic evaluation of three-dimensional wave fields $u$ given by double integrals of the type
\begin{equation}
	\label{eq:Fourier-type-2D_gen}
	u(\bdx) = \iint_{\bdG} F(\bdxi)\exp\{-i r G(\bdxi,\tilde{\bdx})\} \, d\bdxi ,
\end{equation}
where $\bdx \equiv (x_1,x_2,x_3)\in\mathbb{R}^3$ and $\bdxi \equiv (\xi_1,\xi_2)\in\mathbb{C}^2$ represent the physical and spectral variables respectively. The scalar $r$ is defined as $r=|\bdx|$, the unit `observation direction' vector $\tilde{\bdx}$ is given by $\tilde{\bdx}=\bdx/r$, and $d\bdxi$ is understood to be $d\xi_1 \wedge d\xi_2$. A important special case of (\ref{eq:Fourier-type-2D_gen}) are two-dimensional wave fields given by
\begin{equation}
	\label{eq:Fouier2D_gen}
	u(x_1, x_2) = \iint_{\bdG} F(\bdxi)\exp\{-i(x_1 \xi_1+x_2 \xi_2)\} \, d\bdxi.
\end{equation}
In what follows, we will refer to (\ref{eq:Fourier-type-2D_gen}) and (\ref{eq:Fouier2D_gen}) as Fourier-type integrals and Fourier integrals respectively. Most of Part I was dedicated to Fourier integrals, but Fourier-type integrals were also considered in some details in \tmcolor{magenta}{Appendix B}\footnote{Sections, equation numbers, or theorems written in \tmcolor{magenta}{magenta} refer to Part I.}. 
%In particular, we studied the function $G$ that arises in diffraction theory and is given by
%\begin{equation}
%	\label{eq:special-G-diff}
%	G(\bdxi,\tilde{\bdx})=\tilde{x}_1 \xi_1+\tilde{x}_2 \xi_2-|\tilde{x}_3| \sqrt{k^2-\xi_1^2-\xi_2^2},
%\end{equation}
%for some wavenumber $k$. For this choice of $G$, the Fourier-type integral (\ref{eq:Fourier-type-2D_gen}) reduces to the Fourier integral (\ref{eq:Fouier2D_gen}) when $x_3=0$. 

Except for their singular sets, the functions $F$ and $G$ are assumed to be holomorphic functions of $\bdxi$ in a neighbourhood of $\mathbb{R}^2$, and assumed to grow at most algebraically at infinity. The surface of integration $\bdG$ coincides with $\mathbb{R}^2$ almost everywhere, except near the singularity sets of $F$ and $G$ where it is slightly {\it indented}. The {\it indentation} procedure is not straightforward in $\mathbb{C}^2$ and discussed in details in {\cite{Part6A}}, where the process of surface deformation is described by {\it the bridge and arrow notation} (first introduced in \cite{Assier2019c}). The key result of \cite{Part6A} is that the asymptotic behaviour of (\ref{eq:Fourier-type-2D_gen}) as $r\to\infty$ or of (\ref{eq:Fouier2D_gen}) as $x_1^2 + x_2^2\to\infty$ is determined by local integration in the neighbourhoods of several {\it special points} of $F$ or $G$. It is proved that outside the neighbourhoods of such points the surface of integration can be deformed in such a way that the integral is exponentially decaying on it as $x_1^2 + x_2^2\to\infty$ (or $r\to\infty$). The latter is referred to as the {\it locality principle} (first  introduced, for 2D complex analysis, in \cite{Ice2021}). The leading terms of the asymptotic estimations of $u$ can then be found by computing some simple integrals.

In Part II we apply the general results of Part I to the specific example of the problem of diffraction by a quarter-plane. This is a well-known canonical diffraction problem, which was approached by many researchers \cite{satterwhite,Samokish2000,BorovikovPolyhedra,babich95,Assier2012b,Assier2018a,Lyalinov2013, Assier2016}. A good review on the subject can be found in \cite{Assier2018a}. Nevertheless, the task of finding a closed-form solution for the problem of diffraction by a quarter-plane remains open. Innovative two dimensional complex analysis techniques that are developed in {\cite{Assier2019a,Assier2019b,kunz2021diffraction}} seem to be promising in that regard. Of specific interest to the present work are {\cite{Assier2012b}},  {\cite{Lyalinov2013}}  and
{\cite{Shanin2012}} where the authors tried to find the far-field asymptotics of the diffracted field. The geometrical theory of diffraction was used in \cite{Assier2012b}, a  Sommerfeld integral approach was used in \cite{Lyalinov2013}, and ray asymptotics on a sphere with a cut (Smyshlyaev's method) were used in \cite{Shanin2012}. In Part II we recover all the results of these papers, and provide additional formulae for the secondary diffracted field in the plane of the scatterer that, to our knowledge, have never been published before. 
  
%In order to reduce the far-field analysis to the asymptotic estimation of an integral similar to (\ref{eq:Fouier2D_gen}), the so-called {\it analytical continuation formulae} derived in \cite{Assier2018a} are used.

 As is customary, the quarter-plane problem is reduced to a 2D Wiener-Hopf equation involving two unknown spectral functions of two complex variables. As a result, the wavefield can be written as a Fourier-type integral of the form (\ref{eq:Fourier-type-2D_gen}).  The appropriate function $F$ is analytically continued into certain complex domains using the {\it analytical continuation formulae} derived in \cite{Assier2018a}. These formulae provide all the information about the singular set of $F$ needed to asymptotically estimate the wave field using the methods of Part I.

The rest of the article is organised as follows. In section~\ref{sec:sec2formulation}, we provide a mathematical formulation of the problem of diffraction of a plane wave by a quarter-plane with Dirichlet boundary conditions. We indicate that some angles of incidence do not lead to secondary diffracted waves, and some do. We refer to the associated problems as  ``simple'' and ``complicated'', respectively. We introduce the 2D Wiener-Hopf equation in section \ref{sec:2DWH}. We recall some results of \cite{Assier2018a} and write down the {\it analytical continuation formulae} in section \ref{sec:anal-cont-formulae}. Section~\ref{sec:section3} is dedicated to the simple case. We briefly sketch the scheme of the stationary phase method as it was described in Part I (section~\ref{sec:summary-stationary}), we study the singularities of the relevant spectral functions (section~\ref{sec:3.2}) and the indentation of the surface of integration (section~\ref{sec:3.3}), we find the special points (sections~\ref{sec:3.4}--\ref{sec:3.5}) and we build the far-field asymptotics (sections~\ref{sec:3.6}--\ref{sec:3.8}). A similar approach is followed in section~\ref{sec:section4} for the complicated case. We recover all the known asymptotic formulae for the quarter-plane problem, and, in addition, we obtain new formulae for the secondary diffracted waves.

%\begin{itemizedot}
%%  \item Refer to Part 6A {\cite{Part6A}} and summarize the work briefly.   
%%  \item Generic statement about wave diffraction, canonical problems,
%%  quarter-plane as open problem.
%  
%%  \item Talk about {\cite{Assier2012b}}, {\cite{Shanin2012}} and
%%  {\cite{Lyalinov2013}} who all aimed to recover the far-field asymptotics.
%%  Say that in fact we get something new for the secondary diffracted waves.
%  
%  \item Talk about multidimensional complex analysis we applied to the
%  problem. In particular about {\cite{Assier2018a}}, and other double
%  Wiener-Hopf stuff ({\cite{Assier2019a,Assier2019b,kunz2021diffraction}}).
%  Mention {\cite{Radlow1965}}. Make Radlow stuff bigger if we want to talk
%  about it in details in the text (\tmcolor{red}{do we want to do that?}).
%  
%%  \item For simplicity, here we will give the asymptotics formula for the wave
%%  field restricted to the $x_3 = 0$ plane containing the quarter-plane.
%\end{itemizedot}

\section{Formulation, Wiener-Hopf equation, and analytic continuation formulae} \label{sec:section2}
\subsection{Formulation of the diffraction problem. ``Simple'' and ``complicated'' cases} \label{sec:sec2formulation}
We consider the problem of wave diffraction by a quarter-plane ($\tmop{QP}$)
subjected to an incident plane wave and Dirichlet boundary conditions. We make
the time-harmonic hypothesis, with the $e^{- i \omega t}$ convention and time
considerations are henceforth suppressed. The problem reduces to finding the
total wave field $u^t$ satisfying the Helmholtz equation and the Dirichlet
boundary condition
\begin{eqnarray*}
  \Delta u^t + k^2 u^t = 0, & u^t |_{\tmop{QP}} = 0, & u^t = u^{\tmop{in}} +
  u,
\end{eqnarray*}
where $u$ will be referred to as the scattered field. The QP and
the incident wave $u^{\tmop{in}}$ are defined by
\begin{eqnarray*}
  \tmop{QP} = \{ \tmmathbf{x} \in \mathbb{R}^3  , \quad \nobracket x_1 \geqslant
  0, x_2 \geqslant 0, x_3 = 0 \} & \infixand & u^{\tmop{in}} = e^{i (k_1 x_1 +
  k_2 x_2 + k_3 x_3)},
\end{eqnarray*}
where $\tmmathbf{x} = (x_1 , x_2 , x_3)$ and
\begin{eqnarray*}
  k_1 = - k \sin (\theta_0) \cos (\varphi_0), & k_2 = - k \sin (\theta_0) \sin
  (\varphi_0), & k_3 = - \sqrt{k^2 - k_1^2 - k_2^2} = - k \cos (\theta_0),
\end{eqnarray*}
for real spherical incident angles $\theta_0$ and $\varphi_0$.

Instead, one can formulate the problem for the scattered field~$u$.
The latter satisfies the Helmholtz equation and obeys the inhomogeneous  
Dirichlet boundary conditions
\begin{equation}
u(x_1, x_2, 0)|_{\tmop{QP}} = - e^{ i (k_1 x_1 + k_2 x_2) }. 
\label{e:sh001}
\end{equation}
One can see from this formulation that $u(\bdx)$
is symmetric with respect to the plane~$x_3=~0$. For this reason, we will restrict our study to $x_3\geq0$.

%%%%%%%%%%%%%%%%%%%%%%%%%%%%%%
%\begin{figure}[h]
%  \centering{\includegraphics[width=0.5\textwidth]{pictures/figS01.jpg}}
%  \caption{Geometry of the quarter plane problem}
%\label{fig:S01}
%\end{figure}
%%%%%%%%%%%%%%%%%%%%%%%%%%%%%%%

As it is known (see e.g.\ \cite{Assier2018a}), 
the problem formulation should be accompanied by Meixner's conditions at the edges and at the vertex, 
and by a radiation condition at infinity. Meixner's conditions are formulated as the requirement of local integrability 
of the energy-like combination $|\nabla u^t|^2 + |u^t|^2$. These conditions 
prevent the appearance of unphysical sources located at the edges or at the 
vertex.

The radiation condition is more complicated to formulate
and it should be discussed in details.  
Usually, this condition is
formulated in the form of the 
\textit{limiting absorption principle}.  
Assume that the wavenumber parameter 
is represented as 
\begin{equation}
k = k_0 + i \varkappa, 
\label{e:sh002}
\end{equation}
where $k_0$
and $\varkappa$ are the real and the imaginary part of $k$, 
$k_0$ is positive real, and $\varkappa$ is positive but small. 
The value $\varkappa$ corresponds to some absorption in the 
medium.
Fix the value $k_0$ and indicate the dependence of the solution 
on $\varkappa$ by $u^t(\tmmathbf{x} ; \varkappa)$.
The limiting absorption principle states that, for $\varkappa>0$, parts of the field decay exponentially as the distance from the 
vertex tends to infinity.  For real $k = k_0$, one should consider the limit 
\[
u^t(\tmmathbf{x}) =  \lim_{\varkappa \to 0} u^t(\tmmathbf{x}, \varkappa)
\]
of the solution. This scheme is well justified, but sometimes it is difficult to 
find which part of the wave field should decay. Note that the incident wave {\em grows\/}
exponentially as the observation point moves in the direction of incidence.

\begin{figure}[h]
  \centering{\includegraphics[width=0.7\textwidth]{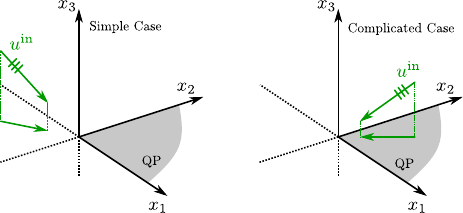}}
  \caption{ Incidence corresponding to the simple case (left) and the complicated case (right)}
\label{fig:S02}
\end{figure}

Let us first consider the {\it simple case\/} by restricting the incident angles to
\begin{equation}
\theta_0 \in (0, \pi / 2)\quad  \mbox{and}  \quad \varphi_0 \in (\pi, 3 \pi / 2),
\label{e:sh003}
\end{equation} 
ensuring that
\begin{eqnarray*}
  {\rm Re}[k_1] > 0 & \infixand & {\rm Re}[k_2] > 0.
\end{eqnarray*}
This case was considered in details in \cite{Assier2018a} and the geometry of the resulting problem is shown in \figurename~\ref{fig:S02}, left.

%%%%%%%%%%%%%%%%%%%%%%%%%%%%%
%%%%%%%%%%%%%%%%%%%%%%%%%%%%%%%

In this simple case, the total far-field is composed of 
the incident plane wave, a reflected plane wave, conical waves emanating from the edges of the quarter plane (the so-called primary diffracted waves), 
and a spherical wave emanating from the vertex (see \figurename~\ref{fig:far-field-simple} for an illustration). There are also penumbral wave fields near the boundaries of the domains occupied by the reflected and the edge diffracted waves.  In particular, we note that in the simple case 
there are no waves that are diffracted by the edges twice (the so-called secondary diffracted waves).

Moreover, for $\varkappa>0$, only the incident wave grows at infinity. The edge diffracted waves are decaying, since the 
projections of the vector $(k_1 , k_2 , k_3)$ onto the edges, i.e.\ $k_1$ and $k_2$, have positive imaginary parts. 
The reflected wave is also decaying. 
Thus, one can formulate the radiation condition as follows: 
the scattered field $u$ should be exponentially decaying at infinity.

Let us now introduce the {\it complicated case} corresponding to
\begin{equation}
\theta_0 \in (0, \pi / 2)\quad  \mbox{and}  \quad \varphi_0 \in (0, \pi / 2),
\label{eq:complicated-angle-restriction}
\end{equation}  
implying that ${\rm Re}[k_1]<0$ and ${\rm Re}[k_2]<0$. The geometry of this case is shown in \figurename~\ref{fig:S02}, right. In addition to the far-field waves present in the simple case, we note the additional presence of secondary diffracted wave in this case (as illustrated in \figurename~\ref{fig:far-field-complicated}). Moreover, in this case, for $\varkappa>0$, the reflected wave and the edge diffracted waves also grow exponentially at infinity. Thus, 
one should extract all these waves (with the penumbral zones) from $u^t$, and require 
that the remaining part of the field is exponentially decaying. This is 
quite complicated and impractical, thus, it would be preferable to have an alternative formulation of the radiation condition. 

A safe way to formulate the radiation condition for the complicated QP problem is considered in \cite{Samokish2000}. 
Instead of a plane wave incidence, one should consider a point source located at the 
point 
\[
\bdx^s=(x_1^{{\rm s}} , x_2^{{\rm s}} , x_3^{{\rm s}}) = R
(\sin(\theta_0) \cos(\vph_0) , \sin(\theta_0) \sin(\vph_0) , \cos(\theta_0)) 
\] 
for some large~$R$ and choose the strength $A$ of the source to be $A = - 4 \pi R\, e^{- i k R}$ to compensate the phase shift and the decay of the incident field.  
Upon denoting by $u^t(\tmmathbf{x} ; \varkappa ,  R)$ the total field resulting from such a source, we require that
for each finite $R$ the field  $u^t(\tmmathbf{x} ; \varkappa , R)$ decays exponentially 
as $|\bdx| \to \infty$. Then $u^t(\tmmathbf{x} ; \varkappa)$ for the plane wave incidence is defined as 
\[
u^t(\tmmathbf{x} ; \varkappa) = \lim_{R \to \infty} u^t(\tmmathbf{x} ; \varkappa ,  R) .
\]
This formulation is mathematically correct\footnote{This is a non-trivial theorem,
since one should prove the existence of the limit.}, but inconvenient 
in practice (the problems with point sources are more complicated than the incident 
plane waves ones).

%Beside the difficulty with the radiation condition, in the complicated case 
%the secondary edge diffracted waves can be emitted.     

In this article, we propose another way to formulate the radiation condition for the 
complicated case. It starts by considering the diffraction problem for the simple case for $\varkappa>0$. 
%Let us solve this problem. 
%or derive any 
%relation for the solution. 
The solution depends on the parameters $k_1$ and $k_2$ that describe the incident wave. For the simple case, these parameters both have a positive real part and a small positive imaginary part. We claim that the solution 
depends analytically on the parameters $k_1$ and $k_2$, and thus it 
%(or any relation for the solution)
remains valid after an analytical continuation
as $k_1$ and $k_2$ are moved along some continuous contours. 
In order to get the complicated case solution,
one should continuously change $k_1$ and $k_2$ until they both have negative real parts and small negative imaginary parts. 

%In more details, we below continue the formulae
%(\ref{eq:third-anal-formula-W}), (\ref{eq:fourth-anal-formula-W})
%from the ``simple'' case to the ``complicated'' case.

The analytical dependency of the solution $u$ on $k_1$ and $k_2$ is a conjecture 
and, formally, one should prove the agreement between the formulation based on the point source described above 
and the analytical continuation. We will not focus on this proof in the present work.
However, one can think of two indirect confirmations that the formulation based on the 
analytical continuation should be correct: 

\begin{itemize}

\item[a.] This is a boundary value problem with boundary data (\ref{e:sh001}) that depends analytically on the parameters $k_1$ and $k_2$ and one would expect that the solution should also depend analytically on these parameters.
%One can consider the problem (\ref{e:sh001})
%for complex parameters  $k_1$ and $k_2$ 
%(rather than for real parameters $\theta_0$, $\vph_0$). 
%The problem is physically admissible, and the solution should be analytical as
%any physical field.  

\item[b.] In what follows, we use this analytic continuation approach to estimate the 
far-field components for the complicated case. The wave components 
found in this way agree with the GTD approach \cite{Assier2012b} and Smyshlyaev's method \cite{Shanin2012}, and no physically prohibited components appear. 

\end{itemize}

Beside the simple and complicated cases, one can study an {\it intermediate}
case, say   $\theta_0 \in (0, \pi / 2)$ and  $\varphi_0 \in (\pi / 2 , \pi)$, for which $\text{Re}[k_1]>0$ and $\text{Re}[k_2]<0$, and only one edge leads to a secondary diffracted wave. 
We will not focus on this, since the methods developed for the complicated case can easily be applied to the intermediate case.  

\subsection{Spectral functions and the 2D Wiener--Hopf equation}  \label{sec:2DWH}
In this section we reproduce briefly the double Fourier transforms manipulations described in
\cite{Assier2018a}. Let us consider the simple case (\ref{e:sh003}) and assume that $\varkappa > 0$. We introduce the double Fourier transform operator $\mathfrak{F}$ defined for any function $\phi$ by 
\begin{equation}
  \mathfrak{F} [\phi (x_1, x_2)] (\tmmathbf{\xi})  = 
  \iint_{\mathbb{R}^2} \phi (x_1 , x_2) e^{i(x_1 \xi_1 + x_2 \xi_2)}\, \mathd x_1  \mathd x_2,
  \qquad  
  \tmmathbf{\xi} \equiv (\xi_1 , \xi_2),
\label{e:sh004}
\end{equation}
and the {\em spectral functions\/} $U(\tmmathbf{\xi})$ and $W(\tmmathbf{\xi})$
defined by 
\begin{align}
U  (\tmmathbf{\xi}) &=\mathfrak{F} [u (x_1 , x_2 , 0^+; \varkappa)] (\tmmathbf{\xi}),
\label{e:sh005}  \\
W  (\tmmathbf{\xi}) &=\mathfrak{F} [\ptl_{x_3}
u (x_1 , x_2 , 0^+; \varkappa)
] (\tmmathbf{\xi}).
\label{e:sh006} 
\end{align}
The dependence of $U$ and $W$ on $\varkappa$ is implied but not indicated. 
Note that the integrals (\ref{e:sh005}) and (\ref{e:sh006}) converge 
due to the limiting absorption principle for the simple case. 
Upon introducing the {\it kernel} 
\begin{equation}
K (\tmmathbf{\xi}) = \frac{1}{\sqrt{k^2 - \xi_1^2 - \xi_2^2}},
\label{e:sh007}
\end{equation}
and assuming that $U$ and $W$ are known, the wave field $u(\tmmathbf{x} ; \varkappa)$
can be reconstructed as 
\begin{equation}
u(\tmmathbf{x} ; \varkappa) = 
\frac{1}{4 \pi^2} \iint_{\mathbb{R}^2}
U(\tmmathbf{\xi})
e^{-i ( \xi_1 x_1 + \xi_2 x_2 - |x_3| / K(\tmmathbf{\xi}))}
\mathd \tmmathbf{\xi}. 
\label{e:sh008}
\end{equation}
The square root in $K$ is chosen to have a positive imaginary part and 
(\ref{e:sh008}) can be interpreted as a plane wave decomposition in the domains $x_3 > 0$
and $x_3 < 0$. 
%Note that the solution $u$ of the problem (\ref{e:sh001}) is symmetric 
%with respect to the symmetry $ x_3 \to - x_3$. 
The differential $\mathd \tmmathbf{\xi}$ denotes simply $\mathd \xi_1 \mathd \xi_2$ 
while $\tmmathbf{\xi}$ belongs to the real plane $\mathbb{R}^2$ and becomes 
$\mathd \tmmathbf{\xi} = \mathd \xi_1  \wedge d\xi_2$
when we start to consider complex integration surfaces. 
The representation (\ref{e:sh008}) enables one to link $U$ and $W$ by the functional equation: 
\begin{equation}
  K (\tmmathbf{\xi}) W (\tmmathbf{\xi})  =  i U (\tmmathbf{\xi}), 
  \label{eq:doubleWHfunctionaleq}
\end{equation} 
and to write down an alternative representation of $u(\tmmathbf{x} , \varkappa)$:
\begin{equation}
u(\tmmathbf{x} ; \varkappa) = 
-\frac{i}{4 \pi^2} \iint_{\mathbb{R}^2}
K(\tmmathbf{\xi})
 W(\tmmathbf{\xi})
e^{-i ( \xi_1 x_1 + \xi_2 x_2 - |x_3| / K(\tmmathbf{\xi}))}
\mathd \tmmathbf{\xi}. 
\label{e:sh009}
\end{equation}
Similarly, the normal derivative of the field on the plane 
$x_3 = 0$ can be written:
\begin{equation}
\ptl_{x_3} u(x_1 , x_2 , 0^+ ; \varkappa) = 
\frac{1}{4 \pi^2} \iint_{\mathbb{R}^2}
 W(\tmmathbf{\xi})
e^{-i ( \xi_1 x_1 + \xi_2 x_2 )}
\mathd \tmmathbf{\xi}.
\label{e:sh010}
\end{equation}

The relation (\ref{eq:doubleWHfunctionaleq}) plays an important role and is actually a 2D Wiener--Hopf functional equation. Indeed, the function $u(x_1 , x_2 , 0^+; \varkappa)$ is known on the 
quarter-plane and is unknown on the remaining 3/4-plane. Conversely, 
due to the symmetry of $u$, the function $\ptl_{x_3} u(x_1 , x_2 , 0^+ ; \varkappa)$
is unknown on the quarter-plane, and is zero 
on the remaining 3/4-plane. Thus, (\ref{eq:doubleWHfunctionaleq}) 
links the Fourier transforms of functions that are unknown on non-overlaping domains whose union is the whole $x_3=0$ plane. 

\subsection{Analytical continuation formulae} \label{sec:anal-cont-formulae}
A Wiener--Hopf formulation of the type (\ref{eq:doubleWHfunctionaleq})
may potentially lead to a solution of the diffraction problem. Unfortunately, up to date, 
no rigorous solution of this 2D Wiener--Hopf problem is known \cite{Assier2019a,kunz2021diffraction}.   
Instead, some useful properties for the Wiener--Hopf problem have 
been derived in \cite{Assier2018a}. They are the so-called {\em analytical 
continuation formulae}\footnote{Note that the quarter-plane is not the only problem for which such formulae can be derived. Similar formulae were also found for the case of a no-contrast right-angled penetrable wedge \cite{Kunz2023}.}. These formulae are integral representations of the unknown function 
$W(\tmmathbf{\xi})$ defining it in certain complex domains. 
We base our further consideration on these representations. 
Here we list the main results of \cite{Assier2018a} that are used in the current paper. 

Before moving further, we need to note that the kernel $K$ admits some
factorisations given by
\[
  K (\xi_1, \xi_2) = \frac{1}{\gamma (\xi_1, \xi_2) \gamma (\xi_1, - \xi_2)}  
  = \frac{1}{\gamma (\xi_2, \xi_1) \gamma
  (\xi_2, - \xi_1)},
\]
where the function $\gamma$ is given by
\begin{equation}
  \gamma (\xi_1, \xi_2)  =  \sqrt{\sqrt{k^2 - \xi_1^2} + \xi_2} .
\label{e:sh010a}
\end{equation}
It is also useful to define the sets $H^{\pm}$, $\hat{H}^{\pm}$ and $P$ as
follows. $H^{+}$ and $H^{-}$ are domains of a single complex variable, defined as the upper and lower half-planes cut along the curves $h^+$ and $h^-$ defined as follows (see \figurename~\ref{fig:Hpm}, left):
\begin{equation}
h^{\pm}= \{\xi\in\mathbb{C} \text{ \ s.t. \ } \xi = \pm \sqrt{k^2 - \tau^2} \text{ \ for \ } \tau \in \mathbb{R}\}. 
\end{equation}
$\hat H^+$ is the upper half-plane, which is not cut along $h^+$, and $\hat H^-$ is the lower half-plane, which is not cut along $h^-$. These domains are all open, i.e.\ the boundary is not included. Denote the contour along $h^-$ as $P$ (see \figurename~\ref{fig:Hpm}, right). The boundary of $H^-$ is therefore
$\ptl H^- = \mathbb{R}\cup P$, and the boundary of $H^+$ is $\ptl H^+ = \mathbb{R}\cup (-P)$.  

%%%%%%%%%%%%%%%%%%%%%%%%%%%%%%%%%%%%%%%%
\begin{figure}[h]
  \centering{\includegraphics[width=0.9\textwidth]{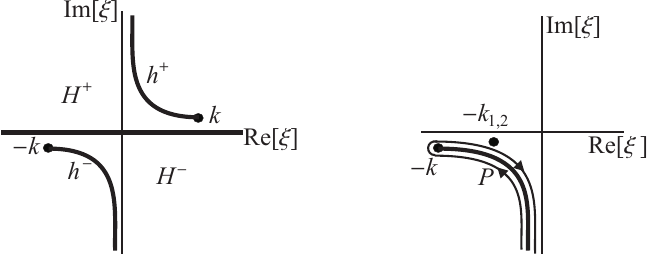}}
  \caption{Domains $H^{\pm}$, left; contour $P$, right }
\label{fig:Hpm}
\end{figure}
%%%%%%%%%%%%%%%%%%%%%%%%%%%%%%%%%%%%%%%%

Consider the formula (\ref{e:sh006}). Initially it defines 
the function $W(\tmmathbf{\xi})$ on the real plane: $\tmmathbf{\xi} \in \mathbb{R}^2$.
However, since $W$ is non-zero only in the first quadrant of the plane, 
the integral converges for any $\tmmathbf{\xi}$ with 
${\rm Im}[\xi_1] > -\epsilon$ and ${\rm Im}[\xi_2] > -\epsilon$, 
where $\epsilon$ can be chosen according to 
\[
0 < \epsilon < \varkappa \cos(\theta_0) \, {\rm min} (\cos(\vph_0) , \sin (\vph_0)). 
\] 
This is based on a rough estimation of decay of the wave field in the plane 
$x_3 = 0$ due to the losses. 
In other words, the integral  (\ref{e:sh006})
converges in $\hat H^+ \times \hat H^+$ and in some 
neighbourhood of the boundary of the real plane\footnote{
In the notation $\hat H^+ \times \hat H^+$, the set before ``$\times$'' is related to the 
$\xi_1$ complex plane, and the set after ``$\times$'' is related to the $\xi_2$ complex plane.}.

In \cite{Assier2018a} the authors pursued the aim to continue $W$ to the remaining parts of the complex space 
of~$\tmmathbf{\xi}$. 
We formulate the main results of \cite{Assier2018a} in the form of two propositions.

\begin{proposition}[Primary continuation]
\label{pr:primary}
The following integral representations for $W$ are valid: 
\begin{align}
	W(\xi_1 , \xi_2) &= \frac{i \gamma(\xi_1, \xi_2) \gamma(\xi_1 , k_2)}{(\xi_1 + k_1)(\xi_2 + k_2)}+
	\frac{\gamma(\xi_1 , \xi_2)}{4\pi^2}
	\int \limits_{-\infty}^{\infty} \mathd \xi_2' 
	\int \limits_{-\infty}^{\infty} \mathd \xi_1'
	\frac{
		\gamma(\xi_1 , - \xi_2') \, K(\xi_1' , \xi_2') \, W(\xi_1' , \xi_2')
	}{ 
		(\xi_1' - \xi_1)(\xi_2' - \xi_2)
	} \label{e:sh011} \\
	W(\xi_1 , \xi_2) &= \frac{i \gamma(\xi_2, \xi_1) \gamma(\xi_2 , k_1)}{(\xi_1 + k_1)(\xi_2 + k_2)}+
	\frac{\gamma(\xi_2 , \xi_1)}{4\pi^2}
	\int \limits_{-\infty}^{\infty} \mathd \xi_1' 
	\int \limits_{-\infty}^{\infty} \mathd \xi_2'
	\frac{
		\gamma(\xi_2 , - \xi_1') \, K(\xi_1' , \xi_2') \, W(\xi_1' , \xi_2')
	}{ 
		(\xi_1' - \xi_1)(\xi_2' - \xi_2)
	}\label{e:sh012}
\end{align}

The formula (\ref{e:sh011}) defines $W$ analytically in the domain $(H^- \setminus \{ - k_1\}) \times \hat H^+$ while the formula (\ref{e:sh012}) defines $W$ analytically in the domain $ \hat H^+ \times (H^- \setminus \{ - k_2\}) $.
\end{proposition}

\begin{proposition}[Secondary continuation]
\label{pr:secondary}
The following integral representations (\ref{eq:third-anal-formula-W}) and (\ref{eq:fourth-anal-formula-W}) for $W$ are valid
in $\hat H^- \times (\hat H^+ \cup H^-)$ and $(\hat H^+ \cup H^-) \times \hat H^-$, respectively:
\begin{align}
W(\xi_1 , \xi_2) &=
\frac{
i \gamma(\xi_1 , \xi_2) \, \gamma(\xi_1 , k_2) \, \gamma(k_2 , k_1)
}{
(\xi_1 + k_1)(\xi_2 + k_2)\, \gamma(k_2 , -\xi_1)
}
+ 
\frac{\gamma(\xi_1 , \xi_2)}{4\pi^2} J_1(\xi_1 , \xi_2),
\label{eq:third-anal-formula-W} \\
J_1 (\xi_1 , \xi_2) &\equiv 
\int \limits_P \mathd \xi_2' 
\int \limits_{-\infty}^{\infty} \mathd \xi_1'
\frac{
\gamma(\xi_1 , -\xi_2') \, K(\xi_1', \xi_2') \, W(\xi_1' , \xi_2')
}{
(\xi_1' - \xi_1)(\xi_2' - \xi_2)
},
\label{e:sh013} \\
W(\xi_1 , \xi_2) &=
\frac{
i \gamma(\xi_2 , \xi_1) \, \gamma(\xi_2 , k_1) \, \gamma(k_1 , k_2)
}{
(\xi_1 + k_1)(\xi_2 + k_2)\, \gamma(k_1 , -\xi_2)
}
+ 
\frac{\gamma(\xi_2 , \xi_1)}{4\pi^2} J_2(\xi_1 , \xi_2),
\label{eq:fourth-anal-formula-W} \\
J_2 (\xi_1 , \xi_2) &\equiv 
\int \limits_P \mathd \xi_1' 
\int \limits_{-\infty}^{\infty} \mathd \xi_2'
\frac{
\gamma(\xi_2 , -\xi_1') \, K(\xi_1', \xi_2') \, W(\xi_1' , \xi_2')
}{
(\xi_1' - \xi_1)(\xi_2' - \xi_2)
}.
\label{e:sh014}
\end{align}
\end{proposition}

Note that to obtain the formulae of Proposition~\ref{pr:secondary}, we need to use Proposition~\ref{pr:primary}, 
since the integrals $J_1$ and $J_2$ require that $W$ is defined on $\mathbb{R} \times P$
and $P \times \mathbb{R}$. 

The formulae (\ref{eq:third-anal-formula-W}) and (\ref{eq:fourth-anal-formula-W}) will be extremely important for the rest of the article, so let us analyse them briefly. They each have an integral 
term and a non-integral term on the right-hand side (RHS). The non-integral terms are easy to deal with, since they 
are products of elementary functions. They are defined everywhere as analytic functions with branch and polar singularities. Due to their simplicity, it is possible to extract their local behaviour near those singularities. 

Conversely, the integral terms contain the unknown function $W(\tmmathbf{\xi})$,
so we can only make general conclusions about them. However, below,  we show that most of the far-field asymptotic wave components of $u$
(all components, except the spherical wave diffracted by the vertex of the QP) arise solely due to the 
non-integral terms of  (\ref{eq:third-anal-formula-W}) and (\ref{eq:fourth-anal-formula-W}).

Let us start by studying the function $J_1(\tmmathbf{\xi})$ defined by (\ref{e:sh013}).
The variables $\xi_1$ and $\xi_2$ can take values in the domains whose boundaries 
are the contours of integration, i.e.\ $\xi_1$ belongs to $\hat H^-$, and 
$\xi_2$ belongs to $\hat H^+ \cup H^-$ (this is the whole plane cut along~$h^-$).
The integral (\ref{e:sh013}) defines a two-sheeted function on this domain, having a branch set 
at $\xi_1 = -k$. This branch set is a complex line with real dimension 2 and comes from the inner square root of the factor $\gamma(\xi_1 , - \xi_2')$. The two sheets of $J_1$ will be referred to as the physical and unphysical sheets. The resulting physical sheet of $W$ is the one used in the formulae (\ref{e:sh009}) and (\ref{e:sh010}). It corresponds to the following choice of the square root of $\sqrt{k^2 - \xi_1^2}$ in $\gamma$:
when $\xi_1$ is real, the values of the square root should be close to positive real or to positive imaginary 
as $\varkappa \to 0$. Note that $J_1$ is regular at all points of the domain 
$\hat H^- \times (\hat H^+ \cup H^-)$ such that $\xi_1 \ne - k$ (on both sheets).

%Being restricted onto the real plane, $J_1$ has the {\em physical sheet\/}, i.e.\ the values 
%participating in the integral 
%(\ref{e:sh010}), and the unphysical sheet. The physical sheet 
%corresponds to the following choice of the square root of $\sqrt{k^2 - \xi_1^2}$ in $\gamma$:
%the values of the square root should be close to positive real or to positive imaginary 
%for $\varkappa \to 0$.  
%The branch set has real dimension~2, i.e.\ it is a complex line.  
%Note that the function $J_1$ is regular at all points of the domain 
%$\hat H^- \times (\hat H^+ \cup H^-)$ having $\xi_1 \ne - k$ (on both sheets).  
  
A similar analysis can be made for the function $J_2$. It is a two-sheeted function in 
the domain 
$ (\hat H^+ \cup H^-) \times \hat H^-$ with a branch set at $\xi_2 = - k$. 

A sketch of the two continuations provided by 
(\ref{eq:third-anal-formula-W}) and (\ref{eq:fourth-anal-formula-W})
is shown in \figurename~\ref{fig:S03}.  
The sketch is made in the coordinates $({\rm Im}[\xi_1] , {\rm Im}[\xi_2])$. The domain~I is the domain of definition of $W$ through (\ref{e:sh006}), the domain~II
corresponds to (\ref{eq:third-anal-formula-W}), and the domain~III to
(\ref{eq:fourth-anal-formula-W}).
   
%%%%%%%%%%%%%%%%%%%%%%%%%%%%%%%%%%%%%%%%
\begin{figure}[h]
  \centering{\includegraphics[width=0.5\textwidth]{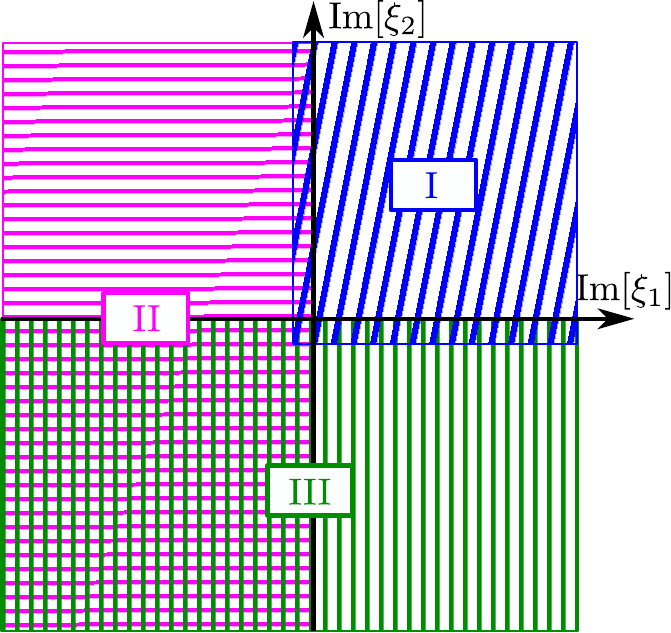}}
  \caption{Continuation of $W$ into the complex domain}
\label{fig:S03}
\end{figure}
%%%%%%%%%%%%%%%%%%%%%%%%%%%%%%%%%%%%%%%%
 
One can see that the domain ${\rm Im}[\xi_1] < 0$, ${\rm Im}[\xi_2] < 0$
is covered by both formulae 
(\ref{eq:third-anal-formula-W}) and (\ref{eq:fourth-anal-formula-W}). 
This is a useful feature. Consider a neighbourhood of some point 
$\tmmathbf{\xi}^\star = (\xi^\star_1 , \xi^\star_2)$
with $\xi^\star_2 \in P$, ${\rm Im}[\xi^\star_1] < 0$. The formula (\ref{eq:third-anal-formula-W}) 
provides a continuation into this neighbourhood, but the neighbourhood is 
cut in two isolated parts by the surface $\xi_2 \in h^-$.
Thus, it is impossible to say whether $\tmmathbf{\xi}^\star$ is a regular point 
of~$W$. Let us assume that $\xi^\star_1 \notin h^-$.
Then one can use the formula (\ref{eq:fourth-anal-formula-W}). This formula provides 
two sheets of continuation of $J_2$ (or $W$) into two samples (on two sheets) of 
the neighbourhood of~$\tmmathbf{\xi^\star}$, and one can study singularities in this neighbourhood. 
%{\color{red} What is the point of this last paragraph?}

One can consider also the points $\tmmathbf{\xi^\star}$ such that 
$\xi_1^\star \in P$ and $\xi_2^\star \in P$. Such points are problematic both for (\ref{eq:third-anal-formula-W}) 
and for (\ref{eq:fourth-anal-formula-W}).
The consideration of these points are as follows. The set $P \times P$ has real 
dimension~2, and {\em it is not an analytical set\/}, i.e.\ it cannot be defined as 
$g (\tmmathbf{\xi}) = 0$ for some holomorphic function~$g$. Thus, by a well-known theorem
of multidimensional complex analysis (Hartogs' theorem, see \cite{Shabat2} p226)
this set, or any 2D neighbourhood on it, cannot be a singularity of~$W$. 
%{\color{red} Same here, where are we going to use this fact?}

\section{Far-field asymptotics in the simple case}
\label{sec:section3}
%%%%%%%%%%%%%%%%%%%%%%%%%%%%%%%%%%%%%%%%%%%%%%

Our aim is to find the far-field asymptotics for the real wavenumber quarter-plane problem in the \textit{simple case}. It means that the incident angles will be restricted by (\ref{e:sh003}), implying that $k_1>0$ and $k_2>0$. To do this, we will apply the machinery of \cite{Part6A}, with the help of the analytical continuation formulae of Section \ref{sec:anal-cont-formulae}.

\subsection{The scheme of the stationary phase method} \label{sec:summary-stationary}

Our aim is to estimate 
the integrals (\ref{e:sh009}) and (\ref{e:sh010}) in the far field based on our knowledge of the singularities of~$W$. Recall that $r$ and the unit observation direction vector $ \tmmathbf{\tilde x}$ are defined as $r=|\bdx|$ and $\tmmathbf{\tilde x} = (\tilde x_1 , \tilde x_2 , \tilde x_3) = \tmmathbf{x} / r$.
%\begin{equation}
%r = \sqrt{x_1^2 + x_2^2 + x_3^2}, 
%\qquad 
%\tmmathbf{\tilde x} = (\tilde x_1 , \tilde x_2 , \tilde x_3) = \tmmathbf{x} / r. 
%\label{e:sh022}
%\end{equation}
Our aim is to get an asymptotic estimation of the field as $ \tmmathbf{\tilde x}$ remains constant, 
and~$r \to \infty$. We consider the limit $\varkappa \to 0$ and only study 
the non-vanishing  terms, i.e.\ the components of the field
that are not exponential decaying as $r \to \infty$.

Upon introducing the functions
\begin{equation}
F (\tmmathbf{\xi}) = K(\tmmathbf{\xi})\, W(\tmmathbf{\xi}), \quad \text{and} \quad	G(\tmmathbf{\xi} ; \tmmathbf{\tilde x}) = \tilde x_1 \xi_1 + \tilde x_2 \xi_2 - |\tilde x_3| / K(\xi_1 , \xi_2), 
	\label{e:sh023}
\end{equation}
the integral (\ref{e:sh009}) can be rewritten as 
\begin{equation}
	u(\tmmathbf{x}; \varkappa) = - \frac{i}{4 \pi^2} \iint_{\mathbb{R}^2}
	F (\tmmathbf{\xi})
	\exp \{ - i r G(\tmmathbf{\xi} ; \tmmathbf{\tilde x}) \} \, \mathd \tmmathbf{\xi},
	\label{e:sh024}
\end{equation}
which is a Fourier-type integral of the form (\ref{eq:Fourier-type-2D_gen}). The integral (\ref{e:sh010})
	is a simpler Fourier integral of the form (\ref{eq:Fouier2D_gen}), and (\ref{e:sh024}) also reduces to such a Fourier integral when $x_3=0$, since $rG(\bdxi;\tilde{\bdx})$ reduces to $\xi_1 x_1+\xi_2 x_2$ in that case. We now list some results about stationary phase method of \cite{Part6A}.

%and (\ref{e:sh024}) becomes a double Fourier integral of the form (\ref{eq:Fouier2D_gen}).
%Introduce the function 
%\begin{equation}
%G(\tmmathbf{\xi} ; \tmmathbf{\tilde x}) = \tilde x_1 \xi_1 + \tilde x_2 \xi_2 - |\tilde x_3| / K(\xi_1 , \xi_2). 
%\label{e:sh023}
%\end{equation}
%{\color{red}For (\ref{e:sh010}) we assume that $\tilde x_3 = 0$ in this formula.} Using $G$,
%the integral (\ref{e:sh009}) can be rewritten as 
%\begin{equation}
%u(\tmmathbf{x}; \varkappa) = - \frac{i}{4 \pi^2} \iint_{\mathbb{R}^2}
%F (\tmmathbf{\xi})
%\exp \{ - i r G(\tmmathbf{\xi} ; \tmmathbf{\tilde x}) \} \, \mathd \tmmathbf{\xi},
%\qquad 
%F (\tmmathbf{\xi}) = K(\tmmathbf{\xi})\, W(\tmmathbf{\xi}),
%\label{e:sh024}
%\end{equation}
%which is a standard form of a phase integral with a large parameter. {\color{red}The integral (\ref{e:sh010})
%can be reduced to a similar form. }

%One should distinguish the Fourier integrals (such as (\ref{e:sh009}) with $x_3 = 0$ and (\ref{e:sh010})) and the plane wave decompositions (such as (\ref{e:sh009})
%with $x_3 > 0$). 
%
%{\color{red} The statements for the Fourier integrals are as follows.}
%
%Statements that are valid for both Fourier-type and Fourier integrals are as follows:
%{\color{red}STOPPED HERE}

\noindent $\bullet$ The method is applicable to Fourier (\ref{eq:Fouier2D_gen}) and Fourier-type (\ref{eq:Fourier-type-2D_gen}) integrals provided that the functions $F$ and $G(\cdot;\tilde{\bdx})$ are holomorphic in some neighbourhood of the real plane except for their singularities $\sigma_j$
(polar or branch sets). 
These singularities are 2D analytic sets in $\mathbb{C}^2$. They should have the \textit{real property}:
their real traces $\sigma'_j=\mathbb{R}^2 \cap \sigma_j$ should be curves in the real plane (rather than points).

\smallskip

\noindent $\bullet$ The non-vanishing field components can be obtained by estimating the integrals (\ref{eq:Fourier-type-2D_gen}) or (\ref{eq:Fouier2D_gen})
near real {\em special points\/} $\tmmathbf{\xi}^\star$ of two types: saddles on singularities (SoS) and crossings of  singularities. 
The SoS are points $\tmmathbf{\xi}^\star$ at which the vector 
\begin{equation}
(\tilde x_1 , \tilde x_2) \text{ for (\ref{eq:Fouier2D_gen})} \quad \text{ or } \quad 	
\nabla G \equiv \left( 
\frac{\ptl G}{\ptl \xi_1} , 
\frac{\ptl G}{\ptl \xi_2}
\right) \text{ for (\ref{eq:Fourier-type-2D_gen})}
\label{e:sh025}
\end{equation}  
is orthogonal to some singularity trace $\sigma'_j$.
Note that with the definition (\ref{e:sh023}), $\nabla G \to (\tilde x_1 , \tilde x_2)$ as $\tilde x_3 \to 0$.
%This statement is the locality principle for the Fourier integrals.   
  
  \smallskip

\noindent $\bullet$ Each such special point may provide or not provide a non-vanishing term. 
This can be established by studying the mutual orientation of the vector $\nabla G$ (for (\ref{eq:Fourier-type-2D_gen})) or $(\tilde x_1 , \tilde x_2)$ (for (\ref{eq:Fouier2D_gen})) at the special point, and the indentation of the integration surface with respect to the singularities (see \tmcolor{magenta}{Section~4} of \cite{Part6A} for details). 
The special points that provide non-vanishing terms are referred to as 
{\em active}.

\smallskip

\noindent $\bullet$ For Fourier integrals (\ref{eq:Fouier2D_gen}), some crossings of singularities do not provide non-vanishing terms for any $\tilde{\bdx}$.
This is the case for {\em additive crossings\/}. A crossing $\tmmathbf{\xi}^\star$ between two real traces of some singularities 
$\sigma_1'$ and $\sigma_2'$ is said to be additive if the function
$F$ can be written near $\tmmathbf{\xi}^\star$ as 
\begin{equation}
F(\tmmathbf{\xi}) = F_1(\tmmathbf{\xi}) + F_2(\tmmathbf{\xi}),
\label{e:sh026}
\end{equation}  
where $F_1$ is singular only at $\sigma_1$, and $F_2$ is singular only at~$\sigma_2$.  The same is true for Fourier-type integrals (\ref{eq:Fourier-type-2D_gen}) provided that $\bdxi^\star$ is not also a crossing for the singularities of $G$.

\smallskip

\noindent $\bullet$  Tangential touch between two real traces do not provide a non-vanishing term unless they are simultaneously a SoS.

\smallskip

The following statements are only relevant for Fourier-type integrals (\ref{eq:Fourier-type-2D_gen}):

\smallskip

\noindent $\bullet$ For the specific choice of $G$ given in (\ref{e:sh024}), only the special points located inside the circle $\sigma_c'=\{\bdxi\in\mathbb{R}^2, \ \xi_1^2+\xi_2^2=k^2\}$ can provide 
non-vanishing terms of the field. This is because the 
term $|x_3| / K$ in the combination $r G$ provides an exponential decay outside $\sigma_c'$.

\smallskip

\noindent $\bullet$ Beside the SoS and crossings of singularities, there is one more type   
of special points: 2D saddle points $\bdxi^\star$ defined such that 
\begin{equation}
\nabla G(\bdxi^\star) = (0,0),
\label{e:sh027}
\end{equation}
which can also give non-vanishing far-field components.

\smallskip

In what follows, we will consider the singularities of the integrands of 
(\ref{e:sh009}) and (\ref{e:sh010}), find the special points, and apply the corresponding asymptotic
formulae obtained in \cite{Part6A}.

%%%%%%%%%%%%%%%%%%%%%%%%%%%%%%%%%%%%%%%%%%%%%%%%%%%%%%%%%%%%%%%%%%%%%% 
 
\subsection{Analysis of the singularities of $W$ near the real plane}\label{sec:3.2}

The formulae of analytical continuation (\ref{eq:third-anal-formula-W}) and
(\ref{eq:fourth-anal-formula-W}) are all we need from \cite{Assier2018a}. They 
provide the required information about the singularities of $W$ near the $(\xi_1 , \xi_2)$ real plane. 
We are specifically interested in the singularities whose intersection with this real plane
has dimension~1 for $\varkappa = 0$, and who belong to the physical sheet defined above.
As we show in \cite{Part6A}, the study of these singularities enables one to estimate non-vanishing components of the 
wave field $u$, say, by using the representation (\ref{e:sh009}).

We denote the singularities by the symbol $\sigma$ with some indexes. Corresponding real traces of 
singularities will be denoted by $\sigma'$ with indexes. Namely, $\sigma'$ 
is $\sigma \cap \mathbb{R}^2$ taken for $\varkappa = 0$.

We will now list these singularities for $W$ and comment on some of their important features. 

\smallskip

\noindent $\bullet$ The set 
\begin{equation}
\sigma_{p_1} = \{ 
\tmmathbf{\xi} \in \mathbb{C}^2, \quad \xi_1 = -  k_1  
\}
\label{e:sh015a}
\end{equation}
is a polar set of $W$. This follows from the non-integral term of (\ref{eq:third-anal-formula-W}).
The integral term of (\ref{eq:third-anal-formula-W}) is not singular at~$\sigma_{p_1}$. 
One can define the residue at $\sigma_{p_1}$ as a function ${\rm res}[W, \sigma_{p_1}](\xi_2)$ by the 
Leray method \cite{Shabat2}. Since $\sigma_{p_1}$ is of the form $\xi_1 = {\rm const}$, this can be found as 
a usual 1D residue in the $\xi_1$ plane for some fixed~$\xi_2$. Namely, 
\begin{equation}
{\rm res}[W, \sigma_{p_1}](\xi_2) = 
\frac{i \gamma (k_1, \xi_2) \gamma(k_1 , k_2)}{\xi_2 + k_2} \cdot
\label{e:sh015}
\end{equation}   
The same residue can be obtained from the non-integral term of (\ref{eq:fourth-anal-formula-W}). Note that the residue has a pole at $\xi_2 = - k_2$ and a branch point at 
$\xi_2 = - \sqrt{k^2 - k_1^2}$. The real trace of $\sigma_{p_1}$ is \\ 
$
\sigma'_{p_1} = 
\{ 
\tmmathbf{\xi} \in \mathbb{R}^2, \quad \xi_1 = - \lim_{\varkappa \to 0} k_1  
\}$.

\smallskip

\noindent $\bullet$ Similarly, the set 
\begin{equation}
\sigma_{p_2} = \{ 
\tmmathbf{\xi} \in \mathbb{C}^2, \quad \xi_2 = - k_2  
\}
\label{e:sh016a}
\end{equation}
is also a polar set of $W$ with residue 
\begin{equation}
{\rm res}[W, \sigma_{p_2}](\xi_1) = 
\frac{i \gamma (k_2, \xi_1) \gamma(k_2 , k_1)}{\xi_1 + k_1} \cdot
\label{e:sh016}
\end{equation}   
This residue has a pole at $\xi_1 = - k_1$ and a branch point at 
$\xi_1 = - \sqrt{k^2 - k_2^2}$. The real trace of $\sigma_{p_2}$ is 
$
\sigma'_{p_2} = 
\{ 
\tmmathbf{\xi} \in \mathbb{R}^2, \quad \xi_2 = - \lim_{\varkappa \to 0}  k_2
\}$.

\smallskip

\noindent $\bullet$ The set 
\begin{equation}
\sigma_{b_1} = \{ 
\tmmathbf{\xi} \in \mathbb{C}^2, \quad \xi_1 = - k 
\}
\label{e:sh017a}
\end{equation}
is a branch set of $W$ of order~2. It is a branch set for the factors
$\gamma(\xi_1 , \xi_2)$, $\gamma(\xi_1 , k_2)$, $\gamma(\xi_1 , \xi_2)$
in (\ref{eq:third-anal-formula-W}). Near $\sigma_{b_1}$, the behaviour of $W$ is given by
\begin{equation}
W(\tmmathbf{\xi}) = C_1 (\xi_2) + C_2(\xi_2) \sqrt{\xi_1 + k} + O(\xi_1 + k). 
\label{e:sh017}
\end{equation}    
The real trace of $\sigma_{b_1}$ is 
$
\sigma'_{b_1} = 
\{ 
\tmmathbf{\xi} \in \mathbb{R}^2, \quad \xi_1 = - k_0.
\}$.

Note that the complex line $\xi_1 = k$ is not a branch set of~$W$. A subset of this line belongs to 
zone~I in \figurename~\ref{fig:S03}, and $W$ should be regular there. 
The rest of this set belongs to zone~III, thus it is described by formula 
(\ref{eq:fourth-anal-formula-W}). The latter has no singularity at $\xi_1 = k$.

\smallskip

\noindent $\bullet$ Similarly, the complex line  
\begin{equation}
\sigma_{b_2} = \{ 
\tmmathbf{\xi} \in \mathbb{C}^2, \quad \xi_2 = - k 
\}
\label{e:sh018}
\end{equation}
is a branch set of order~2 for~$W$. Its real trace is 
$
\sigma'_{b_2} = 
\{ 
\tmmathbf{\xi} \in \mathbb{R}^2, \quad \xi_2 = - k_0
\}$.
The function $W$ is regular on the complex line $\xi_2 = k$. 

\smallskip

\noindent $\bullet$ The analysis of the set 
\begin{equation}
\sigma_c = \{ 
\tmmathbf{\xi} \in \mathbb{C}^2, \quad \xi_1^2 + \xi_2^2 = k^2 
\}
\label{e:sh019}
\end{equation}
is more complicated. Its real trace $\sigma'_c=\{\bdxi\in\mathbb{R}^2, \quad \xi_1^2+\xi_2^2=k_0^2\}$ is the circle
of radius~$k_0$. As we will now see, $W$ is only singular on a portion of this circle. 

To start with, it follows directly from the definition (\ref{e:sh006}) and the properties of quarter-range Fourier transforms, that the function $W$ is regular on $\sigma_c \cap (\hat H^+ \times \hat H^+)$. 

Let us now consider the domain $\mathbb{A} = H^- \times (\hat H^+ \cup \hat H^-)$ and rewrite (\ref{eq:third-anal-formula-W}) as
\begin{equation}
W(\tmmathbf{\xi}) = 
\gamma(\xi_1 , \xi_2) \left( 
\frac{i \gamma(\xi_1 , k_2) \gamma(k_2 , k_1)}{(\xi_1+ k_1)(\xi_2 + k_2) \gamma(k_2 , - \xi_1)}
+ 
\frac{J_1(\xi_1 , \xi_2)}{4 \pi^2}
\right). 
\label{e:sh020}
\end{equation}
The terms within parentheses do not have singularities in $\mathbb{A} \cap \sigma_c$. Thus, the behaviour of $W$ on $\sigma_c$ is determined by that of the factor $\gamma(\xi_1 , \xi_2)$. Being made of usual functions, $\gamma(\xi_1 , \xi_2)$ is easy to analyse, but its Riemann manifold has a non-trivial structure. Indeed, this function is branching at $\xi_1 = \pm k$, and on $\sigma_c$, but not everywhere: it is possible for 
$\sqrt{k^2 - \xi_1^2} + \xi_2$ not to be zero on~$\sigma_c$.
Focusing on the real trace, one can see that, on the physical
sheet, the fragment of $\sigma'_c$ having ${\rm Re}[\xi_2] < 0$ 
bears branching, while the points of $\sigma'_c$ with ${\rm Re}[\xi_2] > 0$
are regular points of~$W$. 

A similar consideration of (\ref{eq:fourth-anal-formula-W}) leads to the conclusion that
the part of $\sigma'_c$ bearing branching of $W$ is the one with ${\rm Re}[\xi_1] < 0$. 
We can therefore conclude that the part of $\sigma_c'$ where $W$ is singular is given by 
\begin{equation}
\sigma_{cc}' = 
\{ 
\tmmathbf{\xi} \in \mathbb{R}^2, \quad 
\tmmathbf{\xi} \in \sigma_c' , \,\, \xi_1 < 0 , \,\, \xi_2 < 0
\}.
\label{e:sh021}
\end{equation}

Moreover, since the branching of $W$ is resulting from either $\gamma(\xi_1 , \xi_2)$
or $\gamma(\xi_2 , \xi_1)$, we can conclude that, near 
$\sigma'_{cc}$, the function $W$ behaves like $\phi(\tmmathbf{\xi})/K(\bdxi)$, 
for some function $\phi(\tmmathbf{\xi})$ regular on~$\sigma'_{cc}$.
Thus, if $\tmmathbf{\xi}$ bypasses $\sigma'_{cc}$, $W(\tmmathbf{\xi})$ 
changes to $-W(\tmmathbf{\xi})$.

A similar reasoning leads to the fact that the product $F(\bdxi)=K(\tmmathbf{\xi}) W(\tmmathbf{\xi})$ arising in 
(\ref{e:sh009}) is singular on $\sigma'_c \setminus \sigma'_{cc}$ and regular on $\sigma'_{cc}$. 
 
\smallskip

According to the formulae of analytic continuation, this list contains all  the
real traces of singularities on the physical sheet. They are sketched in Figure \ref{fig:field-simple-case}.

%%%%%%%%%%%%%%%%%%%%%%%%%%%%%%%%%%%%%%%%%%%%%%%%%%%%%%%
\subsection{Indentation of the integration surface}\label{sec:3.3}

As $\varkappa \to 0$, the singularities listed above hit the real plane at their real traces. As we demonstrated above, and as depicted in \figurename~\ref{fig:field-simple-case}, the real traces for $F=K W$ and $W$ are made up of 5 components given by
$\sigma'_{p_1}$,  $\sigma'_{p_2}$, 
$\sigma'_{b_1}$,  
$\sigma'_{b_2}$, and  
$\sigma'_c$. 
Thus, one cannot use $\mathbb{R}^2$ as the integration surface in (\ref{e:sh009}) or (\ref{e:sh010}). However for $\varkappa>0$, using the 2D Cauchy theorem \cite{Shabat2}, it is possible to slightly deform the integration surface
in $\mathbb{C}^2$ without changing the value of the integral, as long as no singularities are hit in the process. This is done in a way that when taking the limit $\varkappa\to0$, the singularities do not hit the integration surface anymore. We call this deformation an {\it indentation}  of the integration surface around the singularities. The resulting surface of integration (denoted $\bdG$) coincides with $\mathbb{R}^2$
everywhere except in some neighbourhood of the real traces of the singularities.  

For each singularity, there are only two possible types of indentation. This issue is not simple, 
and has been discussed in details in \cite{Part6A}. The choice of indentation is given by the {\it bridge and arrow} 
notations. The correct choice for the problem at hand, dictated by he limiting process $\varkappa \to 0$,  is shown in \figurename~\ref{fig:field-simple-case}. Applying the technique described in \tmcolor{magenta}{Section 3.4} of {\cite{Part6A}},
it is reasonably straightforward to find the bridge and arrow configuration of
$\sigma_{b_{1, 2}}'$ and $\sigma_{p_{1, 2}}'$, while the bridge configuration
of $\sigma_c'$ is determined by the tangential touch compatibility with
$\sigma_{b_{1, 2}}'$.

%%%%%%%%%%%%%%%%%%%%%%%%%%%%%%%%%%%%%%%%%%%%%%%%%%%%%%%%%

\subsection{Special points for the Fourier integrals (\ref{e:sh009})$|_{x_3= 0}$ and (\ref{e:sh010}) }\label{sec:3.4}

%As we demonstrated above, the real traces of singularities of $W K$ and $W$ are made
%up of 5 irreducible singularities given by
%\[
%\sigma'_{p_1}, 
%\quad 
%\sigma'_{p_2}, 
%\quad 
%\sigma'_{b_1}, 
%\quad 
%\sigma'_{b_2}, 
%\quad 
%\sigma'_c. 
%\]
%They 
%are depicted on \figurename~\ref{fig:field-simple-case}. Note that only a subpart of the circle
%$\sigma_c'$ is actually singular, as was discussed above. Using \tmcolor{magenta}{Section 3.4} of {\cite{Part6A}},
%it is reasonably straightforward to find the bridge and arrow configuration of
%$\sigma_{b_{1, 2}}'$ and $\sigma_{p_{1, 2}}'$, while the bridge configuration
%of $\sigma_c'$ is determined by the tangential touch compatibility with
%$\sigma_{b_{1, 2}}'$. The obtained configuration is displayed on \figurename~\ref{fig:field-simple-case}.

\begin{figure}[h]
	\centering{\includegraphics[width=0.9\textwidth]{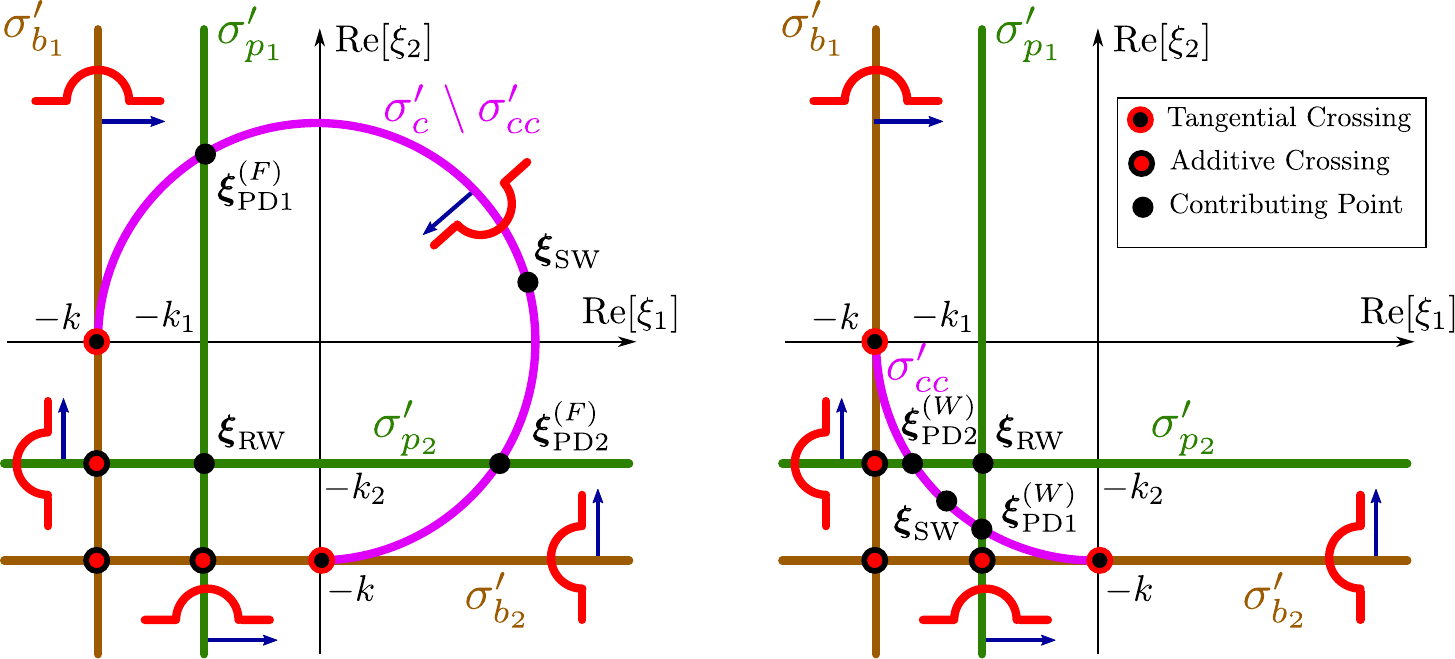}}
	\caption{The real traces and the special points of the functions $F$ (left) and $W$ (right).  
	}
	\label{fig:field-simple-case}
\end{figure}

As discussed in section \ref{sec:summary-stationary}, to study the far-field behaviour of these Fourier integrals, it is enough to look at special points. Those are the SoS and the intersections of real traces as depicted in \figurename~\ref{fig:field-simple-case}. 

We will first focus on the special points that will not contribute to the far-field asymptotics. Let us start with the intersections. The points corresponding to $\sigma'_{b_1}
\cap \sigma'_c$ and $\sigma'_{b_2} \cap \sigma'_c$ intersect tangentially,
and, hence, by \tmcolor{magenta}{Theorem 4.9} of {\cite{Part6A}}, they do not
contribute to the asymptotic expansion of $u$ or $\partial_{x_3}u$. The points
$\sigma'_{p_1} \cap \sigma'_{b_2}$ 
and 
$\sigma'_{p_2} \cap \sigma'_{b_1}$
are transverse crossings with the additive
crossing property. 
This can be proven as follows. Consider the polar line $\sigma_{p_1}$. 
The residue on this line is given by (\ref{e:sh015}). One can see that this residue has no 
branching at the crossing point $\sigma'_{p_1} \cap \sigma'_{b_2}$.  
This means that the pole can be addititively separated from
the branching near the crossing point. 
Thus, the points 
$\sigma'_{p_1} \cap \sigma'_{b_2}$ 
and 
$\sigma'_{p_2} \cap \sigma'_{b_1}$
do not contribute to the asymptotic expansion of $u(x_1, x_2 , 0^+)$ and
$\ptl_{x_3} u(x_1, x_2 , 0^+)$. In Appendix~\ref{app:proof-of-additive-crossing} 
we also show that  that the point
$\sigma'_{b_1} \cap \sigma'_{b_2}$ 
is  an additive crossing
and, therefore, does not contribute to the asymptotic expansions. 

Let us now consider the potential SoS.
In principle, the vector $(\tilde x_1 , \tilde x_2)$ can be orthogonal 
to any of the real traces. However, 
we note that, except for $\sigma'_{c}$, all real traces are straight lines. If $(\tilde x_1 , \tilde x_2)$ is orthogonal to one of these straight traces at one point, it is orthogonal to it everywhere. Such a pathological case 
was excluded from the method developed in \cite{Part6A}, so we 
have to exclude such directions from our consideration. Physically, such directions are likely to belong to the penumbral zones, and mathematically, they would necessitate a non-local approach.

Therefore the relevant special points are the other crossings and the SoS on $\sigma'_c$. For a given direction $(\tilde x_1 , \tilde x_2)$, there are two possible SoS where  $(\tilde x_1 , \tilde x_2)$ is orthogonal to $\sigma_c'$. Obviously, it first needs to be on a singular part of $\sigma'_c$. Moreover for it to count (to be active in the terms of \cite{Part6A}), given the bridge and arrow configuration, $(\tilde x_1 , \tilde x_2)$ attached to this point would need to point towards the origin (see \cite{Part6A} for more detail). An active SoS can hence only be of the form $\tmmathbf{\xi}_{\tmop{SW}}  = - k (\tilde x_1 , \tilde x_2)$. All the relevant special points are therefore
\begin{alignat}{3}
  \tmmathbf{\xi}_{\tmop{RW}} &= (- k_1, - k_2), \quad &&\tmmathbf{\xi}_{\text{PD}1}^{(F)} = (- k_1, k_2'),\quad  &&\tmmathbf{\xi}_{\text{PD}2}^{(F)} = (k_1', -
  k_2),\\
  \tmmathbf{\xi}_{\tmop{SW}} &= - k (\tilde x_1 , \tilde x_2), \quad &&\tmmathbf{\xi}_{\text{PD}1}^{(W)} = (- k_1, - k_2'), \quad &&\tmmathbf{\xi}_{\text{PD}2}^{(W)} = (- k_1',
  - k_2),
\end{alignat}
where
\begin{eqnarray}
  k_2' \equiv \sqrt{k^2 - k_1^2} & \tmop{and} & k_1' \equiv \sqrt{k^2 - k_2^2}
 \label{eq:defofk12prime}
\end{eqnarray}
are both strictly positive. The position of the SoS $\tmmathbf{\xi}_{\tmop{SW}}$ depends solely on the observation direction $(\tilde x_1 , \tilde x_2)$, while the position of all the other special points (crossings) depend solely on the incident angles. Anticipating our findings, the subscripts
$_{\tmop{RW}, \text{PD} 1, \text{PD} 2, \tmop{SW}}$ stand for reflected wave,
primary diffracted wave 1 and 2, and spherical wave.

\subsection{Special points for the Fourier-type integral (\ref{e:sh009})$|_{x_3 > 0}$}\label{sec:3.5}
\begin{figure}[h]
  \centering{\includegraphics[width=0.8\textwidth]{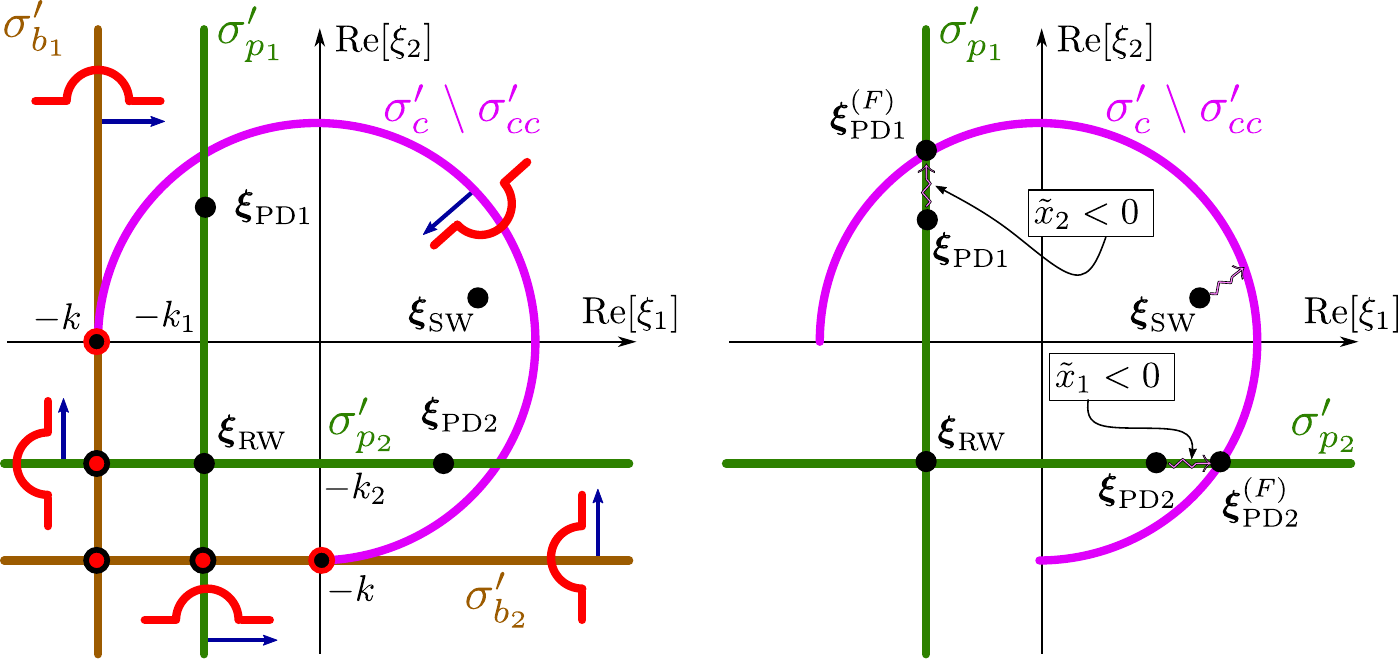}}
  \caption{The real traces of $F$ and special points for (\ref{e:sh009})$|_{x_3 > 0}$ in the simple case (left); behaviour of these points as $\tilde{x}_3\to0^+$ (right) indicated by the symbol $\rightsquigarrow$. }  
\label{fig:field-simple-case_b}
\end{figure}

As discussed in section \ref{sec:summary-stationary} and in more detail in \tmcolor{magenta}{Appendix~B} of \cite{Part6A}, the situation becomes slightly different for $x_3 >0$. In particular, potential special points lying on $\sigma'_c$ stop to provide non-vanishing contributions or become inactive. Moreover, since we have a Fourier-type integral, we also need to consider potential 2D saddle points. For the specific case under investigation, the following is observed. The point $\bdxi_{\text{SW}}$ that was a SoS on $\sigma'_c$ for $x_3=0$ becomes a 2D Saddle and is given by the same formula. However, since $\tilde{x}_3>0$, this point now lies strictly inside the circle $\sigma'_c$. The transverse crossings $\bdxi_{\text{PD}1,2}^{(F)}$ migrate to SoS on $\sigma'_{p_{1,2}}$, again strictly inside $\sigma'_c$. Note that now a SoS is defined as a point where $\nabla G$ is perpendicular to a real trace. Since this vector also depend on $\bdxi$, it is now possible to have an isolated SoS on a straight real trace. The crossing $\bdxi_{\text{RW}}$ remains contributing, without a change in type or location. The special points to consider and their type (displayed in \figurename~\ref{fig:field-simple-case_b}) are hence given by
\begin{align}
	\text{(Crossing)}&: \tmmathbf{\xi}_{\tmop{RW}} = (- k_1, - k_2), \quad \quad \, \text{(2D saddle)}: \tmmathbf{\xi}_{\tmop{SW}} = - k (\tilde x_1 , \tilde x_2), \nonumber \\
	\quad\text{(SoS)}&: \tmmathbf{\xi}_{\tmop{PD1}} =  \left(-k_1,\frac{-\tilde x_2 k_2'}{\sqrt{{\tilde x}_2^2 +{\tilde x}_3^2}}\right), \quad \text{(SoS)}: \tmmathbf{\xi}_{\tmop{PD2}} =  \left(\frac{-\tilde x_1 k_1'}{\sqrt{{\tilde x}_1^2 +{\tilde x}_3^2}},-k_2\right). \label{eq:special-points-3D-simple}
\end{align}

In the next section, to each of these special points, we will associate a far-field wave component by using the results in
\tmcolor{magenta}{Section 5} of
{\cite{Part6A}}. For brevity, we will only do it for the field $u$, that is for the integral (\ref{e:sh009}). However, we have given all the necessary results for the method to also be applied to get the far-field asymptotics for $\partial_{x_3}u(x_1,x_2,0^+)$ using the integral (\ref{e:sh010}).

\subsection{Field estimation near the special points for (\ref{e:sh009})$|_{x_3 > 0}$}\label{sec:3.6}
 
\subsubsection{The reflected wave}
Consider the vicinity of  $\tmmathbf{\xi}^{\star} =\tmmathbf{\xi}_{\tmop{RW}}$,
which is a transverse crossing point for the real trace components $\sigma'_1
= \sigma'_{p_1}$ and $\sigma_2' = \sigma'_{p_2}$. We can approximate $F$  near this point using formula (\ref{eq:third-anal-formula-W}) 
to obtain the local approximations
\begin{eqnarray*}
 F (\tmmathbf{\xi}) \underset{\tmmathbf{\xi} \rightarrow
  \tmmathbf{\xi}_{\tmop{RW}}}{\approx} \frac{i}{ (\xi_1
  + k_1) (\xi_2 + k_2)},\quad
% W (\tmmathbf{\xi}) \underset{\tmmathbf{\xi} \rightarrow
%  \tmmathbf{\xi}_{\tmop{RW}}}{\approx} \frac{i/K(k_1,k_2)}{ (\xi_1
%  + k_1) (\xi_2 + k_2)}.
\end{eqnarray*}
Then, let us suppose that $\nabla G\neq(0,0)$ near $\bdxi^\star$, i.e. there is no stationary point nearby, therefore, 
$G$ can be approximated as follows: 
\begin{equation}
\label{phase_asympt}
G(\bdxi)  \underset{\tmmathbf{\xi} \rightarrow
  \tmmathbf{\xi}_{\tmop{RW}}}{\approx} G(\tmmathbf{\xi}^{\star}) +  (\tmmathbf{\xi}-\tmmathbf{\xi}^\star)\cdot\nabla G(\tmmathbf{\xi}^{\star}).
\end{equation}
%\[
%\frac{\ptl G}{\ptl \xi_1}() = \tilde x_1 + \frac{\tilde x_3\xi_1}{\sqrt{k^2 -\xi_1^2-\xi_2^2}},\quad \frac{\ptl G}{\ptl \xi_2} = \tilde x_2 + \frac{\tilde x_3\xi_2}{\sqrt{k^2 -\xi_1^2-\xi_2^2}}.
%\]
Thus, the associated wave component is given by the following integral: 
\begin{equation}
	\label{uRW_cont}
	u_{\tmop{RW}} (\tmmathbf{x}) = A\iint_{\bdG}\frac{\exp\left\{ -ir(\bdxi-\bdxi^\star)\cdot\nabla G(\bdxi^\star)\right\}}{ (\xi_1 + k_1) (\xi_2 + k_2)}d\tmmathbf{\xi},
\end{equation}
where 
\[
A = \tfrac{1}{4\pi^2}e^{-irG(\tmmathbf{\xi}^\star)}, \ G(\bdxi^\star)=-k_1 \tilde{x}_1-k_2 \tilde{x}_2+k_3 \tilde{x}_3, \ \nabla G(\bdxi^\star)=(\tilde{x}_1+\tfrac{\tilde{x}_3 k_1}{k_3},\tilde{x}_2+\tfrac{\tilde{x}_3 k_2}{k_3}). 
\] 
%\begin{equation}
%\label{uRW_cont}
%u_{\tmop{RW}} (\tmmathbf{x}) = A\iint_{\bdG}\frac{\exp\left\{ ir\frac{\ptl G}{\ptl \xi_1}(\tmmathbf{\xi}^{\star}) \xi_1 + ir\frac{\ptl G}{\ptl \xi_2}(\tmmathbf{\xi}^{\star})\xi_2\right\}}{ (\xi_1 + k_1) (\xi_2 + k_2)}d\tmmathbf{\xi},
%\end{equation} 
%where 
%\[
%A = \frac{1}{4\pi^2}\exp\left\{-ir(G(\tmmathbf{\xi}^\star) - \bdxi^\star\cdot\nabla G(\bdxi^\star))\right\},
%\]

The surface of integration $\bdG$ is chosen to coincide with $\mathbb{R}^2$ everywhere except for the two singular lines $\sigma'_{p_{1,2}}$, where it is slightly indented according to bridge notation shown in \figurename~\ref{fig:field-simple-case_b}. Hence, we can use the results of \tmcolor{magenta}{Section 5.2} and \tmcolor{magenta}{Appendix B} of {\cite{Part6A}} to get:
\begin{equation}
\label{eq:uPWsimple}
  u_{\tmop{RW}} (\tmmathbf{x})  = -\exp\{i(k_1x_1+k_2x_2-k_3x_3)\}\mathcal{H} \left(x_1+\frac{x_3k_1}{k_3}\right)\mathcal{H} \left(x_2+\frac{x_3k_2}{k_3}\right),
\end{equation}
where  $\mathcal{H}$ is the Heaviside function. Note that the Heaviside functions in
(\ref{eq:uPWsimple})  define the region of activity of the crossing $\tmmathbf{\xi}_{\tmop{RW}}$.
One can also notice that (\ref{uRW_cont}) is just a composition of two one-dimensional polar integrals so there are other, simpler, means to evaluate (\ref{uRW_cont}). This is exactly the result for the reflected wave that one would expect from Geometrical Optics considerations.

\subsubsection{Primary diffracted waves} \label{sec:primary-simple-fourier-type}
Consider the vicinity of $\tmmathbf{\xi}^{\star} =\tmmathbf{\xi}_{\text{PD}
1}$, that corresponds to a SOS on  $\sigma_1' =
\sigma_{p_1}'$.  Using (\ref{eq:third-anal-formula-W}) we can approximate $F$ near $\tmmathbf{\xi}^{\star}$ as follows
%{\color{red}\begin{eqnarray*}
%	F (\tmmathbf{\xi}) \underset{\tmmathbf{\xi} \rightarrow
%		\tmmathbf{\xi}_{\tmop{PD1}}}{\approx}  \frac{\sqrt{i(k_2'+k_2)(x_2^2+x_3^2)}}{(-x_2k'_2 + k_2)\sqrt{k_2'}(\sqrt{x_2^2+x_3^2}-x_2)}\times\frac{1}{\xi_1+k_1},
%\end{eqnarray*}}
%{\color{red} I think the above is plain wrong, it should be}
\begin{eqnarray*}
	F (\tmmathbf{\xi}) \underset{\tmmathbf{\xi} \rightarrow
		\tmmathbf{\xi}_{\tmop{PD1}}}{\approx}  \frac{i\sqrt{k_2'+k_2}(\tilde{x}_2^2+\tilde{x}_3^2)^{3/4}}{\sqrt{k_2'}(k_2 \sqrt{\tilde{x}_2^2+\tilde{x}_3^2}-k_2' \tilde{x}_2)\sqrt{\sqrt{\tilde{x}_2^2+\tilde{x}_3^2}+\tilde{x}_2}}\times\frac{1}{\xi_1+k_1}.
\end{eqnarray*}
Since $\bdxi^\star$ is a SOS, i.e.\ $\nabla G \perp \sigma_1'$ at $\bdxi^\star$, following \tmcolor{magenta}{Section 5.1} and \tmcolor{magenta}{Appendix B} of {\cite{Part6A}}, and noting that $\frac{\ptl G}{\ptl\xi_2}(\tmmathbf{\xi}^{\star})=0$,  we only need to consider the following approximation for $G$:
\begin{equation}
\label{eq:phase_expansion}
G(\bdxi)  \underset{\tmmathbf{\xi} \rightarrow
  \tmmathbf{\xi}_{\tmop{PD1}}}{\approx} G(\tmmathbf{\xi}^{\star}) +  (\xi_1 - \xi^\star_1)\frac{\ptl G}{\ptl \xi_1}(\tmmathbf{\xi}^{\star}) + \frac{(\xi_2 - \xi^\star_2)^2}{2}\frac{\ptl^2 G}{\ptl \xi_2^2}(\tmmathbf{\xi}^{\star}).
\end{equation}
Terms in $\mathcal{O}((\xi_1 - \xi^\star_1)^2)$, $\mathcal{O}((\xi_1 - \xi^\star_1)(\xi_2 - \xi^\star_2))$ and higher orders can be neglected.
Using that
\begin{equation*}
	G(\bdxi^\star)=-k_1 \tilde{x}_1-k_2'\sqrt{\tilde{x}_2^2+\tilde{x}_3^2}, \quad \frac{\partial G}{\partial \xi_1}(\bdxi^\star)=\tilde{x}_1-\frac{k_1}{k_2'}\sqrt{\tilde{x}_2^2+\tilde{x}_3^2}, \quad \frac{\partial^2 G}{\partial \xi_1^2}(\bdxi^\star)=\frac{(\tilde{x}_2^2+\tilde{x}_3^2)^{3/2}}{k_2'\tilde{x}_3^2},
\end{equation*}
we obtain the asymptotic component $u_{\text{PD} 1} (\tmmathbf{x})$ associated to the special point $\bdxi_{\text{PD}1}$:
\begin{align}
  u_{\text{PD} 1} (\tmmathbf{x}) &= \tfrac{e^{-i\frac{3\pi}{4}}x_3\sqrt{k'_2 +k_2}\exp\left\{ik_1x_1 + ik_2'\sqrt{x_2^2 + x_3^2}\right\}}{\sqrt{2\pi}(k_2\sqrt{x_2^2+x_3^2}-x_2k_2')\sqrt{\sqrt{x_2^2+x_3^2}+x_2}}\mathcal{H} \left(x_1-\frac{k_1}{k_2'}\sqrt{x_2^2+x_3^2}\right).
  \label{eq:uPD1simple}
\end{align}
A very similar approach with $\bdxi_{\text{PD}2}$ leads to the asymptotic component $u_{\text{PD} 2}(\bdx)$:
\begin{align}
  u_{\text{PD} 2} (\tmmathbf{x}) &= \tfrac{e^{-i\frac{3\pi}{4}}x_3\sqrt{k'_1 +k_1}\exp\left\{ik_2x_2 + ik_1'\sqrt{x_1^2 + x_3^2}\right\}}{\sqrt{2\pi}(k_1\sqrt{x_1^2+x_3^2}-x_1k_1')\sqrt{\sqrt{x_1^2+x_3^2}+x_1}}\mathcal{H} \left(x_2-\frac{k_2}{k_1'}\sqrt{x_1^2+x_3^2}\right), 
  \label{eq:uPD2simple}
\end{align}
which can also be recovered by just swapping the indices 1 and 2 in (\ref{eq:uPD1simple}).

\subsubsection{The spherical wave}
Consider the vicinity of the 2D saddle point  $\tmmathbf{\xi}^{\star} =\tmmathbf{\xi}_{\tmop{SW}}$. Apart from some pathological cases (penumbral zones),  $W(\tmmathbf{\xi})$ is known to be regular at $\tmmathbf{\xi}^{\star}$. Therefore, using the fact that $G(\bdxi^\star)\!=\!-k$, $K(\bdxi^\star)\!=\!(k\tilde{x}_3)^{-1}$, $\det(H_G(\bdxi^\star))\!=\!(k\tilde{x}_3)^{-2}$, where $H_G$ is the Hessian matrix of $G$, the classical 2D saddle point method leads to the asymptotic component $u_{\tmop{SW}}$:
\begin{equation}
\label{eq:uSWsimple}
u_{\tmop{SW}}(\tmmathbf{x}) = -\frac{kW(\tmmathbf{\xi}^{\star})}{2\pi}\frac{e^{ikr}}{kr}.
\end{equation}
Note that, unlike the  asymptotic formulae derived previously, this one is not
completely in closed-form. Indeed, it depends on the quantity $W
(\tmmathbf{\xi}^{\star})$, which is not known {\tmem{a priori}}. Means of
finding its value were discussed at length in {\cite{Smyshlyaev1990,shanin1,Assier2012,Assier2019a}}.

\subsection{On the consistency with the asymptotics of (\ref{e:sh009})$|_{x_3 = 0}$}\label{sec:3.7}
As we mentioned above, the Fourier integral (\ref{e:sh009})$|_{x_3 = 0}$ and the Fourier-type integral (\ref{e:sh009})$|_{x_3 > 0}$ need to be treated separately because the special points of (\ref{e:sh009}) change their type as $x_3\to 0$. However, physically, such limit should not be pathological and we expect the resulting asymptotics to be consistent. As $x_3=0$, we can apply the results of \cite{Part6A}  \textit{ad-hoc} without the need of \tmcolor{magenta}{Appendix~B}.
\subsubsection{The reflected wave}
Neither the type (transverse crossing) nor the position of the special point $\tmmathbf{\xi}^{\star} =\tmmathbf{\xi}_{\tmop{RW}}$ depend on $x_3$, so we do not expect any changes. Indeed, using \tmcolor{magenta}{Section 5.2} of \cite{Part6A}, we find that its contribution is given by
\begin{equation}
u_{\tmop{RW}}(x_1,x_2,0) = -\exp\{i(k_1x_1+k_2x_2)\} \mathcal{H}(x_1)\mathcal{H}(x_2),
\end{equation} 
which is consistent with (\ref{eq:uPWsimple}).

\subsubsection{Primary diffracted waves}
For $x_3=0$, the special point $\tmmathbf{\xi}^{\star} =\tmmathbf{\xi}^{(F)}_{\tmop{PD1}}$ is a transverse crossing, and we can hence apply \tmcolor{magenta}{Section 5.2} of \cite{Part6A}. Using ({\ref{eq:third-anal-formula-W}), $F$ can be shown to behave as:
\begin{equation}
F (\tmmathbf{\xi}) \underset{\tmmathbf{\xi} \rightarrow
  \tmmathbf{\xi}^{(F)}_{\tmop{PD1}}}{\approx} \frac{i\sqrt{2k'_2}}{\sqrt{k'_2+k_2}}\times\frac{1}{(\xi_1+k_1)\sqrt{k^2-\xi_1^2-\xi_2^2}}, 
\end{equation}
leading to the far-field wave component
\begin{equation}
u_{\tmop{PD1}}(x_1,x_2,0) = \frac{e^{-i\frac{3\pi}{4}}\exp\{i(k_1x_1-k'_2x_2)\}}{\sqrt{k'_2+k_2}\sqrt{-\pi x_2}}\mathcal{H}(-x_2)\mathcal{H}(k'_2x_1+k_1x_2).
\label{eq:uPD1-Fourier-x3-zero}
\end{equation}
Considering the limit $\tilde{x}_3\to0$ in Section \ref{sec:primary-simple-fourier-type}, we observe the following. If $\tilde{x}_2<0$, as illustrated in figure \ref{fig:field-simple-case_b} (right) then $\bdxi_{\text{PD}1} \to \bdxi_{\text{PD}1}^{(F)}$ and the contribution (\ref{eq:uPD1simple}) $\to$ (\ref{eq:uPD1-Fourier-x3-zero}). If however $\tilde{x}_2>0$, then $\bdxi_{\text{PD1}}\to\bdxi_{\text{PD1}}^{(W)}$, which is not a special point of (\ref{e:sh009})$|_{x_3 = 0}$, and the contribution tends to zero. This behaviour illustrates the perfect consistency between this section and the results of Section \ref{sec:primary-simple-fourier-type} as $x_3\to0$.
%Because as $x_3\to0$, $x_3/\sqrt{\sqrt{x_2^2+x_3^2}+x_2}\to\sqrt{-2x_2}\mathcal{H}(-x_2)$, and, as indicated in \figurename~\ref{fig:field-simple-case_b} (right) $\bdxi_{\text{PD}1} \to \bdxi_{\text{PD}1}^{(F)}$ if $\tilde{x}_2<0$, the latter is consistent with (\ref{eq:uPD1simple}). If $\tilde{x}_2>0$, $\bdxi_{\text{PD1}}\to\bdxi_{\text{PD1}}^{(W)}$, which is not a special point of (\ref{e:sh009})$|_{x_3 = 0}$.

Contribution from the transverse crossing $\tmmathbf{\xi}^{(F)}_{\tmop{PD2}}$ and the consistency with Section \ref{sec:primary-simple-fourier-type} can be obtained in a similar way, or by swapping the indices 1 and 2 in (\ref{eq:uPD1-Fourier-x3-zero}). 
\subsubsection{The spherical wave}
When $x_3=0$, the special point $\tmmathbf{\xi}^{\star} =\tmmathbf{\xi}_{\tmop{SW}}$ is a SOS. Taking into account that $W$ is regular at $\tmmathbf{\xi}^{\star}$, and using the definition (\ref{e:sh023}, left), $F$ is approximated as follows:
\begin{equation}
F (\tmmathbf{\xi}) \underset{\tmmathbf{\xi} \rightarrow
  \tmmathbf{\xi}_{\tmop{SW}}}{\approx}  W(\tmmathbf{\xi}^{\star})\frac{1}{\sqrt{k^2-\xi_1^2-\xi_2^2}}.
\end{equation}
Then, using the results of
\tmcolor{magenta}{Section 5.1} of {\cite{Part6A}} the resulting far-field component is:
\begin{equation}
u_{\tmop{SW}}(x_1,x_2,0) = -\frac{kW(\tmmathbf{\xi}^{\star})}{2\pi}\frac{e^{ik\sqrt{x_1^2+x_2^2}}}{k\sqrt{x_1^2+x_2^2}}\mathcal{H}(\mathbb{R}^2\backslash\tmop{QP}),
\end{equation}
where $\mathcal{H}(\mathbb{R}^2\backslash \tmop{QP})$ is equal to $1$ if $(x_1,x_2)\in \mathbb{R}^2\backslash \tmop{QP}$ and $0$ otherwise. Remembering that $F$ is regular on $\sigma_{cc}'$, one can see that the latter is indeed consistent with (\ref{eq:uSWsimple}).

\subsection{General asymptotic expansion of $u$ }\label{sec:3.8}

We have now considered all the contributing points and we can therefore
reconstruct the far-field asymptotic approximation of $u$ as follows:
\begin{eqnarray*}
  u (\tmmathbf{x}) & \underset{| \tmmathbf{x} | \rightarrow \infty}{\approx} &
  u_{\tmop{RW}} (\tmmathbf{x}) + u_{\text{PD} 1}  (\tmmathbf{x}) +
  u_{\text{PD} 2} (\tmmathbf{x}) + u_{\tmop{SW}} (\tmmathbf{x}),
\end{eqnarray*}
where the components are given by equations (\ref{eq:uPWsimple}),
 (\ref{eq:uPD1simple}), (\ref{eq:uPD2simple}), (\ref{eq:uSWsimple}). These expansions agree with former results given for
example in {\cite{Assier2012b}}, {\cite{Shanin2012}}, {\cite{Lyalinov2013}},
{\cite{Budaev2005a}}, obtained via the GTD, Sommerfeld integrals and the random
walk method. A schematic illustration of these wave components is given in \figurename~\ref{fig:far-field-simple}.

\begin{figure}[h]
	\centering{
		\includegraphics[height=0.4\textwidth]{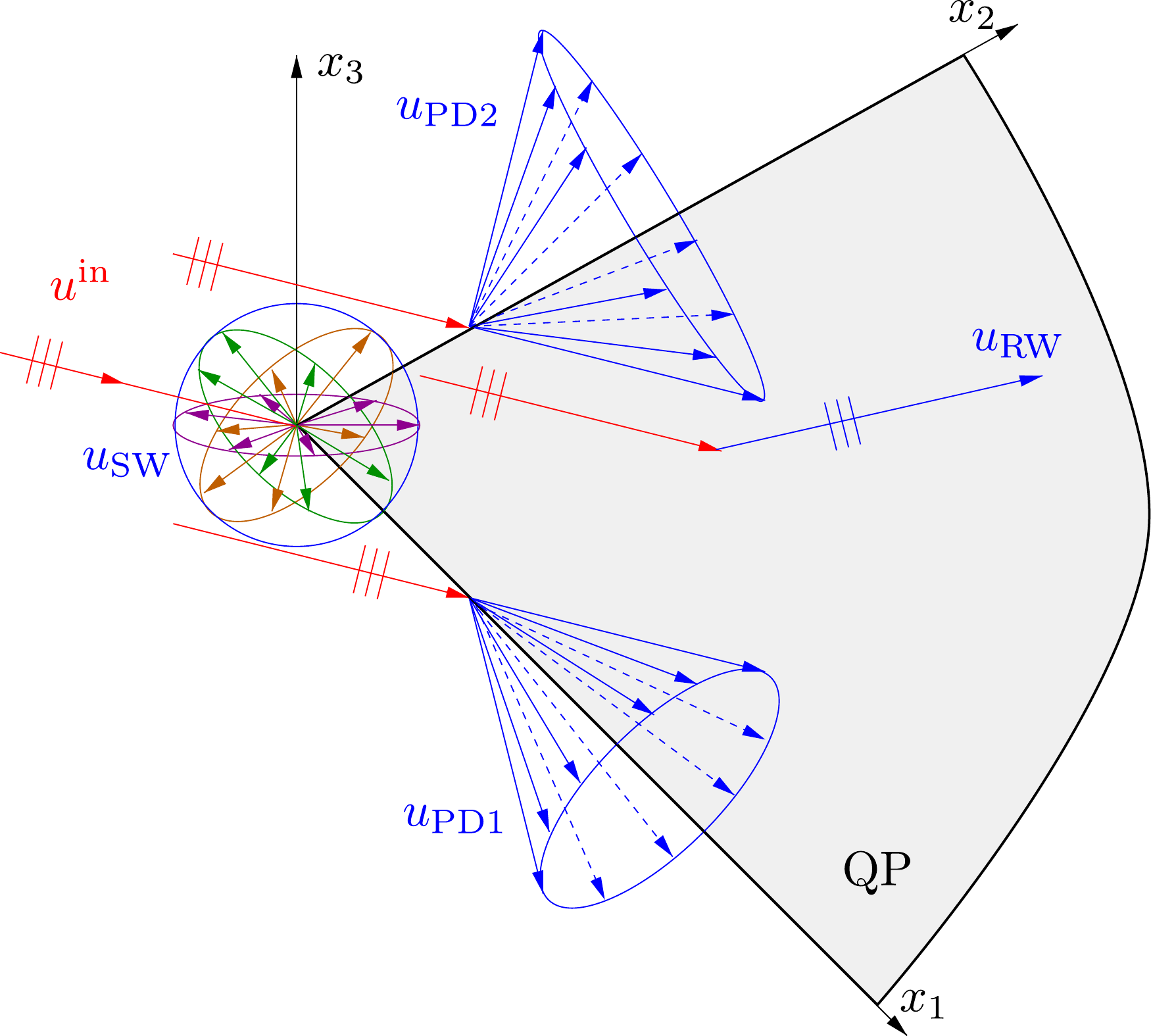}
	}
	\caption{Illustration of the far-field asymptotic components in the simple case (adapted from \cite{Assier2012}).}
	\label{fig:far-field-simple}
\end{figure}

\section{Far-field asymptotics in the complicated case}
\label{sec:section4}
We now consider the \textit{complicated case} for the real wavenumber quarter-plane problem. This means that the incident angles are restricted by (\ref{eq:complicated-angle-restriction}), implying that $k_1<0$ and $k_2<0$. As in the previous section, our aim is to find the resulting far-field asymptotics.
%
%Here we consider the case of incident angles restricted to be such that
%$\theta_0 \in (0, \pi / 2)$ and $\varphi_0 \in (0, \pi / 2)$, ensuring that we
%have
%\begin{eqnarray*}
%  {\rm Re}[k_1] < 0 & \infixand & {\rm Re}[k_2] < 0.
%\end{eqnarray*}
The main reason why this case is labelled \textit{complicated} is
because, on top of the wave fields discussed in the previous section, it is
also known to give rise to two secondary diffracted waves. These are
notoriously difficult to analyse. Moreover, all the known published formulae for these
secondary waves {\cite{Assier2012b,Shanin2012,Lyalinov2013}}, though completely valid for $x_3>0$, have been
obtained in such a way that they blow up on the plane $x_3 = 0$.  In what follows we will show that our method
leads to these same formulae for $x_3>0$, but also provides a proper secondary diffracted wave fields on the $x_3 = 0$
plane, hence providing formulae that, to our knowledge, have never been obtained previously.

\subsection{Special points for the Fourier integrals (\ref{e:sh009})$|_{x_3 = 0}$ and (\ref{e:sh010})}

Assume temporarily that the wavenumber $k$ has a small positive imaginary part~$\varkappa$, and treat ${\rm Re}[k_1]$ and ${\rm Re}[k_2]$ as parameters that can move continuously from a state were
they are both positive, to a state where they are both negative. The free terms of the analytical
continuation formulae exhibit new singular sets when ${\rm Re}[k_{1, 2}]$ change sign.
Indeed, the set  
\begin{equation}
\sigma_{\tmop{sb}_1} = \{ 
\tmmathbf{\xi} \in \mathbb{C}^2, \quad \xi_1 = - k_1'
\},
\end{equation}
where $k_1'$ is defined as in (\ref{eq:defofk12prime}), is a branch line of $W$ of order 2; the subscript $_{\tmop{sb}}$ stands for secondary branch. This is due to the factor $\gamma(\xi_1,k_2)$ in the non-integral term of (\ref{eq:third-anal-formula-W}). 
The function $W$ behaves near $\sigma_{\tmop{sb}_1}$ as 
\begin{equation}
 W(\tmmathbf{\xi}) = \hat C_1(\xi_2) + \hat C_2(\xi_2)\sqrt{\xi_1 + k_1'} + \mathcal{O}(\xi_1 + k_1'). 
\end{equation}
The real trace of $\sigma_{\tmop{sb}_1}$ is 
\begin{equation}
\sigma_{\tmop{sb}_1}' = \{ 
\tmmathbf{\xi} \in \mathbb{R}^2, \quad \xi_1 = - \lim_{\varkappa\to0}k_1' 
\}.
\end{equation}
Note that the factor $\gamma(k_2,-\xi_1)$ in (\ref{eq:third-anal-formula-W}) has a branch point at $k_1'$, but it does belongs to the upper half-plane, where nothing is known about the analyticity properties  of the integral term $J_1$.
Similarly, due to the term $\gamma(\xi_2,k_1)$ in (\ref{eq:fourth-anal-formula-W}), the line 
\begin{equation}
\sigma_{\tmop{sb}_2} = \{ 
\tmmathbf{\xi} \in \mathbb{C}^2, \quad \xi_2 = - k_2' 
\},
\end{equation}
where $k_2'$ is defined in (\ref{eq:defofk12prime}), is a branch line of $W$ of order 2. The real trace of $\sigma_{\tmop{sb}_2}$ is 
\begin{equation}
\sigma_{\tmop{sb}_2}' = \{ 
\tmmathbf{\xi} \in \mathbb{R}^2, \quad \xi_2 = - \lim_{\varkappa\to0}k_2' 
\}.
\end{equation}
One has to be very careful to make sure that the bridge configuration
of the singularities is preserved in this continuous process. This is why, in
a somewhat counter intuitive manner, the bridge configuration of the
singularity $\sigma_{p_{1, 2}}'$ needs to remain the same, even though, for
$\varkappa > 0$, the singularity changes half-plane. This is to make sure that
the singularities do not intersect the surface of integration in the process
of $k_1$ and $k_2$ changing sign. A schematic view of the real traces of the singularity sets of $F$ and $W$,
together with the associated bridge configurations, is given in \figurename~\ref{fig:complicated-case}. 

\begin{figure}[h]
	\centering{\includegraphics[width=0.8\textwidth]{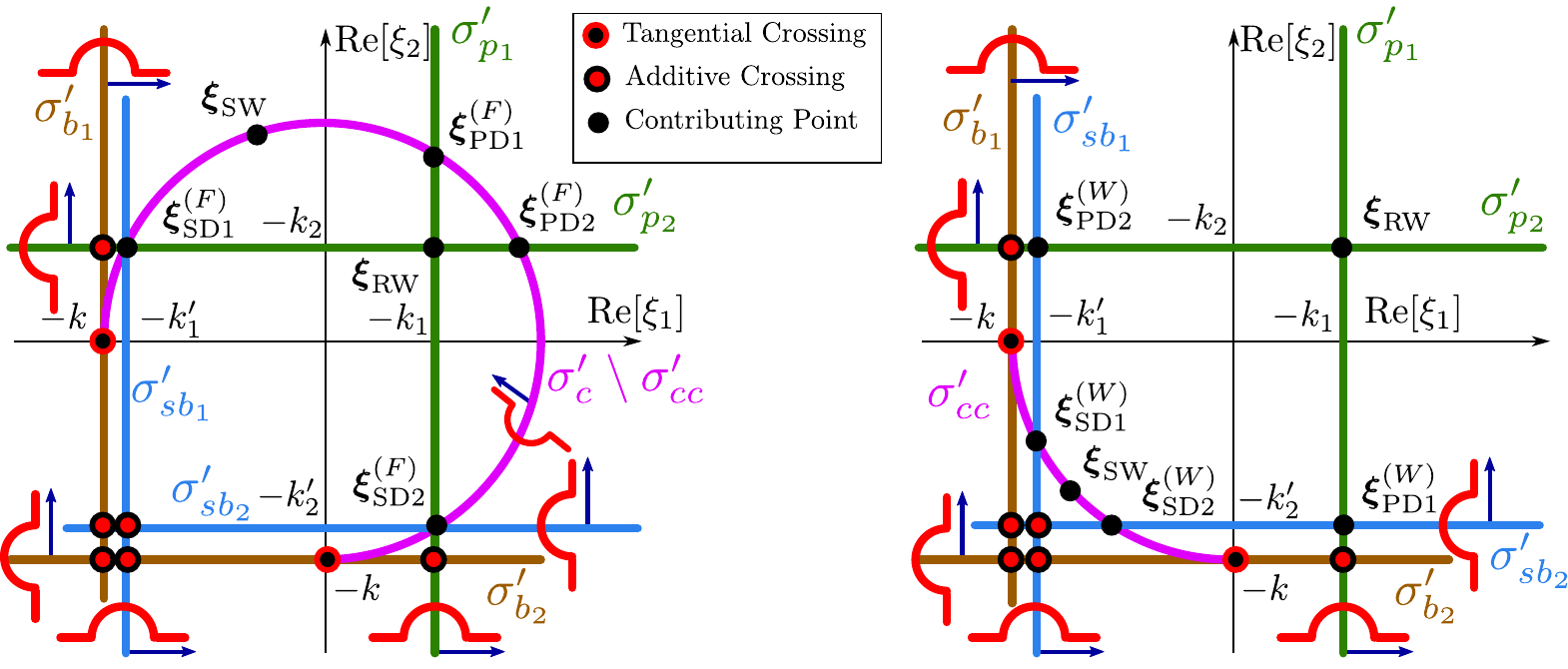}}
	\caption{The real traces and the important points of the spectral function of $F$ (left) and $W$ (right) in the complicated case.}
	\label{fig:complicated-case}
\end{figure}

As for the simple case, and as indicated on this same figure, a few potential special points can be discarded due to either additive or tangential crossing. The additivity of the crossing $\sigma_{\tmop{sb}_1}\cap\sigma_{\tmop{sb}_2}$ follows directly from the structure of the non-integral terms of (\ref{eq:third-anal-formula-W}) or (\ref{eq:fourth-anal-formula-W}), the additivity of crossings $\sigma_{\tmop{b}_2}\cap\sigma_{\tmop{sb}_1}$ and $\sigma_{\tmop{b}_1}\cap\sigma_{\tmop{sb}_2}$ follows from the non-integral terms of (\ref{eq:third-anal-formula-W}) and (\ref{eq:fourth-anal-formula-W}), correspondingly. The special points that need to be considered are:
\begin{alignat*}{3}
	\tmmathbf{\xi}_{\tmop{RW}} &= (- k_1, - k_2), \quad  &&\tmmathbf{\xi}_{\text{PD}
		1}^{(F)} = (- k_1, k_2'),\quad  &&\tmmathbf{\xi}_{\text{PD} 2}^{(F)} = (k_1', -
	k_2),\\
	\tmmathbf{\xi}_{\tmop{SW}}& = -k(\tilde{x}_1, \tilde{x}_2), \quad &&\tmmathbf{\xi}_{\text{PD}
		1}^{(W)} = (- k_1, - k_2'),\quad  &&\tmmathbf{\xi}_{\text{PD} 2}^{(W)} = (- k_1',
	- k_2),\\
	\tmmathbf{\xi}_{\text{SD} 1}^{(F)} &= (- k_1', - k_2), \quad 
	&&\tmmathbf{\xi}_{\text{SD} 2}^{(F)} = (- k_1, - k_2'), \quad  
	&&\tmmathbf{\xi}_{\text{SD} 1}^{(W)} = (- k_1', k_2), \quad
	\tmmathbf{\xi}_{\text{SD} 2}^{(W)} = (k_1, - k_2').
\end{alignat*}

Apart from the SoS $\bdxi_{\text{SW}}$, they are all transverse non-additive crossings. The points $\tmmathbf{\xi}_{\text{SD} 1}^{(F)}$ and $\tmmathbf{\xi}_{\text{SD} 2}^{(F)}$ are remarkable as they correspond to \textit{triple} transverse crossings. 
 
% The points corresponding to $\sigma'_{b_1}
%\cap \sigma'_c$ and $\sigma'_{b_2} \cap \sigma'_c$ still intersect
%tangentially, and, hence, by \tmcolor{magenta}{Theorem 4.9} of
%{\cite{Part6A}}, they do not contribute to the asymptotic expansion of $u$.
%It is shown in Appendix \ref{app:proof-of-additive-crossing} that
% points corresponding to $\sigma_{\tmop{sb}
%1}' \cap \sigma_{\tmop{sb} 2}'$, $\sigma_{\tmop{sb}_1}' \cap \sigma_{b_2'}$
%and $\sigma_{\tmop{sb} 2}' \cap \sigma_{b_1}'$ are transverse crossing with
%the additive crossing property and, therefore, by \tmcolor{magenta}{Theorem
%4.12} of {\cite{Part6A}} they do not contribute to the asymptotic expansion of
%the field.

\subsection{Special points for the Fourier-type integral (\ref{e:sh009})$|_{x_3 > 0}$}
As for the simple case, the nature of the special points change when $x_3>0$. In particular, all the contributing points are now strictly within the circle $\sigma_c'$ and 2D saddles can appear. In this case, the following special points are found:
\begin{alignat*}{3}
	\tmmathbf{\xi}_{\tmop{RW}} &= (- k_1, - k_2), \quad 
	&&\tmmathbf{\xi}_{\text{PD} 1} = \left(- k_1, \frac{-\tilde{x}_2 k'_2}{\sqrt{\tilde{x}_2^2+\tilde{x}_3^2}}\right), \quad 
	&&\tmmathbf{\xi}_{\text{PD} 2} = \left(\frac{-\tilde{x}_1 k'_1}{\sqrt{\tilde{x}_1^2+\tilde{x}_3^2}}, - k_2\right),\\
	\tmmathbf{\xi}_{\tmop{SW}} &= -k(\tilde{x}_1, \tilde{x}_2),\quad  
	&&\tmmathbf{\xi}_{\text{SD} 1} = \left(- k_1', \frac{k_2\tilde{x}_2}{\sqrt{\tilde{x}_2^2+\tilde{x}_3^2}}\right), \quad
	&&\tmmathbf{\xi}_{\text{SD} 2} = \left(\frac{k_1\tilde{x}_1}{\sqrt{\tilde{x}_1^2+\tilde{x}_3^2}}, - k_2'\right).
\end{alignat*}
Apart from $\tmmathbf{\xi}_{\tmop{RW}}$ that remains a transverse crossing and $\tmmathbf{\xi}_{\tmop{SW}}$ that is now a 2D saddle point, all the other special points listed above are now SoS. They are illustrated in \figurename~\ref{fig:complicated-case-x3neq0} (left).
%Special points for  (\ref{e:sh009})$|_{x_3 > 0}$ are shown in \figurename~\ref{fig:complicated-case-x3neq0}.  We denoted  SOS on $\sigma_{\tmop{sb}_{1}}$ and $\sigma_{\tmop{sb}_{2}}$ as  $\tmmathbf{\xi}_{\tmop{SD1}}$ and $\tmmathbf{\xi}_{\tmop{SD2}}$, correspondingly, since they provide asymptotic contributions corresponding to the secondary diffracted waves. Note that for small $x_3$ and $x_1>0$ SOS $\tmmathbf{\xi}_{\tmop{PD2}}$ is close to the triple crossing $\tmmathbf{\xi}^{(F)}_{\tmop{SD1}}$ and merges with  $\tmmathbf{\xi}_{\tmop{SD1}}$ in the limit. A similar process is happening near $\tmmathbf{\xi}^{(F)}_{\tmop{SD2}}$. 
\begin{figure}[h]
  \centering{\includegraphics[width=0.8\textwidth]{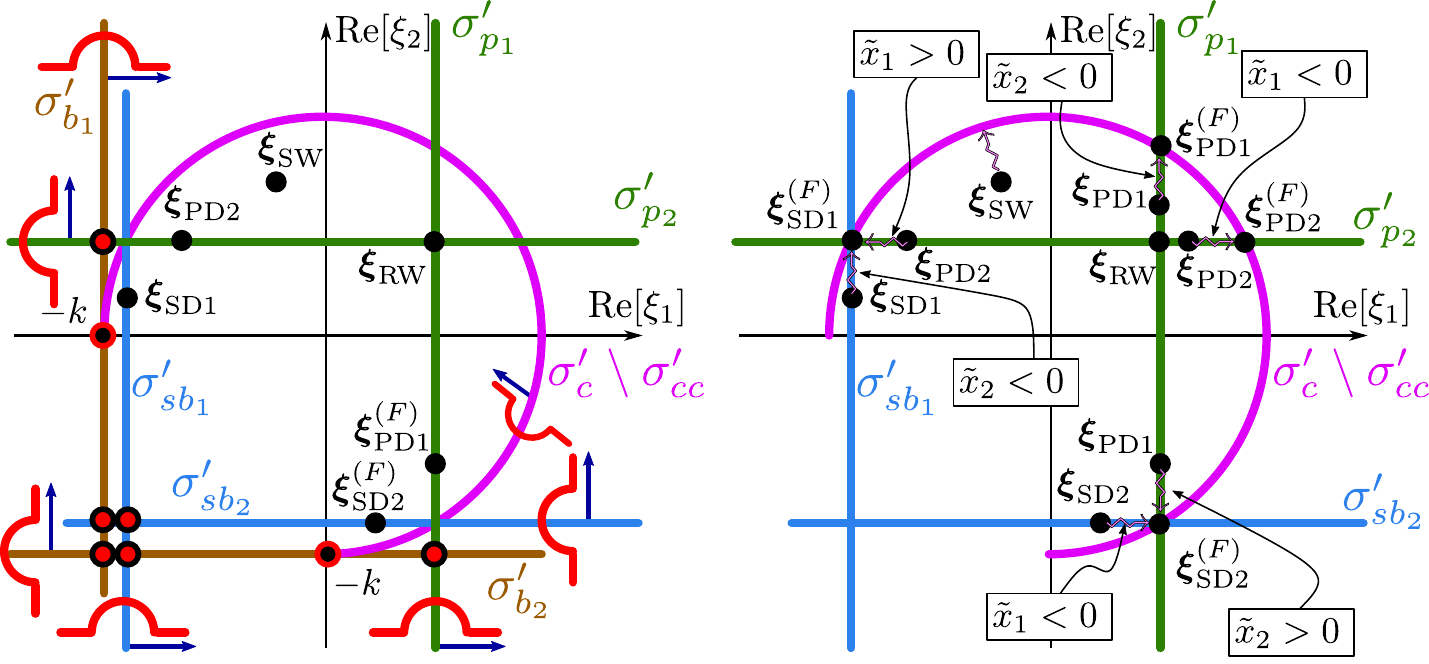}}
  \caption{The real traces of $F$ and the special points of (\ref{e:sh009})$|_{x_3 > 0}$  in the complicated case (left); behaviour of these points as $\tilde{x}_3\to0^+$ (right) indicated by the symbol $\rightsquigarrow$.}
  \label{fig:complicated-case-x3neq0}
\end{figure}

Moreover, as illustrated in \figurename~\ref{fig:complicated-case-x3neq0} (right), we note that
\begin{alignat}{2}
\tmmathbf{\xi}_{\text{PD} 1} &\underset{\tilde{x}_3 \rightarrow
	0}{\longrightarrow} \left\{ \begin{array}{ccc}
	\tmmathbf{\xi}_{\text{PD} 1}^{(F)} & \text{if} & \tilde{x}_2 < 0\\
	\tmmathbf{\xi}_{\text{SD} 2}^{(F)} & \text{if} & \tilde{x}_2 > 0
\end{array} \right., \quad && 	\tmmathbf{\xi}_{\text{PD} 2} \underset{\tilde{x}_3 \rightarrow
0}{\longrightarrow} \left\{ \begin{array}{ccc}
\tmmathbf{\xi}_{\text{PD} 2}^{(F)} & \text{if} & \tilde{x}_1 < 0\\
\tmmathbf{\xi}_{\text{SD} 1}^{(F)} & \text{if} & \tilde{x}_1 > 0
\end{array} \right.,  \\
\tmmathbf{\xi}_{\text{SD} 1} &\underset{\tilde{x}_3 \rightarrow
	0}{\longrightarrow} \left\{ \begin{array}{ccc}
	\tmmathbf{\xi}_{\text{SD} 1}^{(F)} & \text{if} & \tilde{x}_2 < 0\\
	\tmmathbf{\xi}_{\text{SD} 1}^{(W)} & \text{if} & \tilde{x}_2 > 0
\end{array} \right., \quad && \tmmathbf{\xi}_{\text{SD} 2} \underset{\tilde{x}_3 \rightarrow
0}{\longrightarrow} \left\{ \begin{array}{ccc}
\tmmathbf{\xi}_{\text{SD} 2}^{(F)} & \text{if} & \tilde{x}_1 < 0\\
\tmmathbf{\xi}_{\text{SD} 2}^{(W)} & \text{if} & \tilde{x}_1 > 0
\end{array} \right.. \label{eq:limitsx3ofSDicomplicated}
\end{alignat}
%Coordinates of special points are as follows:
%\begin{alignat*}{3}
%\tmmathbf{\xi}_{\tmop{RW}} &= (- k_1, - k_2), \quad 
%&&\tmmathbf{\xi}_{\text{PD} 1} = \left(- k_1, \frac{-\tilde{x}_2 k'_2}{\sqrt{\tilde{x}_2^2+\tilde{x}_3^2}}\right), \quad 
%&&\tmmathbf{\xi}_{\text{PD} 2} = \left(\frac{-\tilde{x}_1 k'_1}{\sqrt{\tilde{x}_1^2+\tilde{x}_3^2}}, - k_2\right),\\
%\tmmathbf{\xi}_{\tmop{SW}} &= -k(\tilde{x}_1, \tilde{x}_2),\quad  
%&&\tmmathbf{\xi}_{\text{SD} 1} = \left(- k_1', \frac{k_2\tilde{x}_2}{\sqrt{\tilde{x}_2^2+\tilde{x}_3^2}}\right), \quad
%&&\tmmathbf{\xi}_{\text{SD} 2} = \left(\frac{k_1\tilde{x}_1}{\sqrt{\tilde{x}_1^2+\tilde{x}_3^2}}, - k_2'\right).
%\end{alignat*}

\subsection{Reflected, primary diffracted and spherical waves}
The points $\bdxi_{\text{RW}}$, $\bdxi_{\text{SW}}$ and $\bdxi_{\text{PD}1,2}$ are given by the same formulae as in the simple case (see (\ref{eq:special-points-3D-simple})) and are of the same type. Therefore the whole procedure of obtaining the
asymptotic components for the corresponding waves remains unchanged and can be
carried out {\tmem{ad-hoc}} in the exact same manner as in the simple case
of Section \ref{sec:section3}, leading to the exact same asymptotic formulae
(\ref{eq:uPWsimple}), (\ref{eq:uPD1simple}),
(\ref{eq:uPD2simple}) and (\ref{eq:uSWsimple}). The same conclusions also hold for these waves when it comes to the consistency between the $x_3=0$ and the $x_3>0$. The only difference is that, as $\tilde{x}_{2,1}>0$, $\bdxi_{\text{PD}1,2} \to \bdxi_{\text{SD}2,1}^{(F)}$ as $\tilde{x}_3\to0$, which this time is a special point (this was not the case in the simple case). However, the asymptotic contribution of $\bdxi_{\text{PD}1,2}$ still tends to zero in that case. More will be said below about the reason for this seeming contradiction.

\subsection{The secondary diffracted waves for (\ref{e:sh009})$|_{x_3 > 0}$}
Consider the vicinity of  $\tmmathbf{\xi}^{\star} =\tmmathbf{\xi}_{\text{SD}
1}$, a SOS on $\sigma_1' = \sigma_{\tmop{sb}_1}'$. Using (\ref{eq:third-anal-formula-W}) we can approximate $F$ near $\tmmathbf{\xi}^{\star}$ as follows:
\begin{equation}
F(\tmmathbf{\xi}) \approx \frac{i\sqrt{k_1'+k_1}(\tilde{x}_2^2+\tilde{x}_3^2)^{3/4}}{\sqrt{2}(k_1'-k_1)k_2^2(\tilde{x}_2+\sqrt{\tilde{x}_2^2+\tilde{x}_3^2})^{3/2}}\times \sqrt{\xi_1 + k_1'}.
\end{equation}
Then, using the expansion (\ref{eq:phase_expansion}) with
\begin{alignat*}{3}
	G(\bdxi^\star)=-k_1'\tilde{x}_1+k_2\sqrt{\tilde{x}_2^2+\tilde{x}_3^2}, \quad && \tfrac{\partial G}{\partial \xi_1}(\bdxi^\star)=\tilde{x}_1+\tfrac{k_1'}{k_2}\sqrt{\tilde{x}_2^2+\tilde{x}_3^2}, \quad && \tfrac{\partial^2G}{\partial \xi_2^2}(\bdxi^\star)=-\tfrac{ \left(\tilde{x}_2^2+\tilde{x}_3^2\right)^{3/2} }{ k_2\tilde{x}_3^2 },
\end{alignat*}	
	 together with the results of
\tmcolor{magenta}{Section 5.1} and \tmcolor{magenta}{Appendix B} of {\cite{Part6A}}, we obtain:
\begin{align}
 \label{eq:usd1formula_x3neq0}
  u_{\text{SD} 1} (\tmmathbf{x}) &=  \tfrac{-x_3 \sqrt{k'_1+k_1} \exp\left\{i\left(k'_1x_1-k_2\sqrt{x_2^2+x_3^2}\right)\right\}}{4\pi(k'_1-k_1)\left(x_2 + \sqrt{x_2^2+x_3^2}\right)^{3/2}\left(-k_2x_1-k'_1\sqrt{x_2^2+x_3^2}\right)^{3/2}}\mathcal{H} \left(-k_2x_1 - k'_1\sqrt{x_2^2+x_3^2}\right).
\end{align} 

This formula agrees with previously published work \cite{Assier2012b,Shanin2012,Lyalinov2013}.
It is interesting to see what happens to this expression as $x_3\to0$. Note first that due to the Heaviside function in (\ref{eq:usd1formula_x3neq0}), this wave component is zero if $x_1<0$. We hence only need to consider the case $x_1>0$.

If we also have $x_2>0$, the expression (\ref{eq:usd1formula_x3neq0}) tends to zero as $x_3\to0$. This is hardly surprising since, in that case, according to (\ref{eq:limitsx3ofSDicomplicated}), $\bdxi_{\text{SD}1}\to\bdxi_{\text{SD}1}^{(W)}$, which is not a special point of (\ref{e:sh009})$|_{x_3=0}$. It is also a confirmation that $u_{\text{SD}1}$ satisfies the Dirichlet boundary conditions on QP.

If instead we have $x_2<0$, the expression blows up in the limit $x_3\to0$, which is more surprising and somewhat unphysical. This is something that was not realised in \cite{Assier2012b,Shanin2012,Lyalinov2013}. Within the framework of the present article, the mathematical reason for this behaviour is due to the fact that as $x_1>0$, $x_2<0$ and $x_3\to0$, two special SOS points ($\tmmathbf{\xi}_{\text{SD} 1}$ and $\tmmathbf{\xi}_{\text{PD} 2}$) merge into the triple crossing $\tmmathbf{\xi}_{\text{SD} 1}^{(F)}$. This arbitrary proximity between these two SoS means that their respective contribution cannot be computed independently as we did and one should treat them together to obtain an accurate picture. We can also give a physical interpretation to this phenomena. To obtain the secondary diffracted wave emanating from the $x_1$ edge, one can approximate, locally, the primary diffracted wave coming from the $x_2$ edge as a plane wave hitting the $x_1$ edge with a grazing incidence along the plate. As can be understood by considering a simpler 2D (half-plane) problem, the region $x_3=0$ outside the plate corresponds to a penumbral zone for this problem. In the next section we will address this issue by deriving a formula valid on the plane $x_3=0$.

%Note that the latter blows up as $x_3\to 0$. From a mathematical point of view it happens due to  merging of two SOS ($\tmmathbf{\xi}_{\text{SD} 1}$ and $\tmmathbf{\xi}_{\text{PD} 2}$) that leads to a triple crossing  $\tmmathbf{\xi}_{\text{SD} 1}^{(F)}$. We study the vicinity of this point in the next section. This phenomenon also has  clear physical interpretation. Indeed, as $x_3 \to 0$ we are getting close to the half-shadow zone where the wavefronts of the primary and secondary diffracted fields are touching.

 A similar reasoning about $\tmmathbf{\xi}^{\star} =\tmmathbf{\xi}_{\text{SD}
 	1}$ , leads to $u_{\text{SD} 2}$, with similar remarks as $x_3\to0$:

\begin{align}
 \label{eq:usd2formula_x3neq0}
 \!\!\! u_{\text{SD} 2} (\tmmathbf{x}) =  \tfrac{-x_3\sqrt{k'_2+k_2} \exp\left\{i\left(k'_2x_2-k_1\sqrt{x_1^2+x_3^2}\right)\right\}}{4\pi(k'_2-k_2)\left(x_1 + \sqrt{x_1^2+x_3^2}\right)^{3/2}\left(-k_1x_2-k'_2\sqrt{x_1^2+x_3^2}\right)^{3/2}}\mathcal{H} \left(-k_1x_2 - k'_2\sqrt{x_1^2+x_3^2}\right).
\end{align}

\subsection{The secondary diffracted waves for (\ref{e:sh009})$|_{x_3 = 0}$}
Consider the vicinity of
$\tmmathbf{\xi}^{\star} =\tmmathbf{\xi}_{\text{SD} 1}^{(F)} = (- k_1', -
k_2)=(\xi_1^\star,\xi_2^\star)$, that corresponds to the transverse intersection of the three singular
traces $\sigma_1' = \sigma_{\tmop{sb}_1}'$, $\sigma_2' = \sigma_c'$ and
$\sigma_3'=\sigma'_{p_2}$.

The case of a triple transverse crossing was not considered in {\cite{Part6A}}
and hence the asymptotics  cannot be
obtained by a direct application of this paper. We will however adapt our method to obtain the asymptotics in this specific case.  
Using  (\ref{eq:third-anal-formula-W}), the leading singular behaviour of $F$ near $\tmmathbf{\xi}^{\star}$ can be found to be:
\begin{equation}
 F (\tmmathbf{\xi})  \underset{\tmmathbf{\xi} \rightarrow
  \tmmathbf{\xi}_{\text{SD} 1}^{(F)}}\approx  \frac{i\sqrt{k_1' + k_1}}{ (- k_1' + k_1)}  \times \frac{(\xi_1 + k_1')^{1 /2}}{(\xi_2 + k_2) (k^2 - \xi_1^2 - \xi_2^2)^{1 / 2}} \cdot
\end{equation}
Upon introducing the local change
of variables $\Psi : (\xi_1, \xi_2) \rightarrow (\rho, \psi)$ defined by
\begin{eqnarray*}
  \xi_1-\xi_1^\star = \rho \cos (\psi) & \tmop{and} & \xi_2-\xi_2^\star = \rho \sin
  (\psi),
\end{eqnarray*}
there is an equivalence between $\bdxi\to\bdxi^\star$ and $\rho\to0$, and $\mathd \tmmathbf{\xi} \leftrightarrow \rho \mathd \rho \wedge \mathd
\psi$. Because $k_2^2+(k_1')^2=k^2$, we can introduce the angle $\vartheta_1$ defined by 
\begin{eqnarray}
	k_2=-k\cos(\vartheta_1) & \text{ and } & k_1'=k\sin(\vartheta_1).
	\label{eq:deftheta1}
\end{eqnarray}
We can hence write
\begin{align*}{2}
  %K (\tmmathbf{\xi})
  k^2 - \xi_1^2 - \xi_2^2 & \underset{\tmmathbf{\xi} \rightarrow
  \tmmathbf{\xi}^{\star}}{=} 2 k_1' \rho \cos (\psi) + 2 k_2 \rho \sin
  (\psi) + \mathcal{O} (\rho^2)=2k\sin(\vartheta_1-\psi)+\mathcal{O}(\rho^2).
\end{align*}
Moreover, using the polar coordinates $(x_1,x_2)=(r \cos (\varphi), r \sin (\varphi))$, we find that 
\begin{align*}
	x_1 \xi_1+x_2 \xi_2&=	x_1 \xi_1^\star+x_2 \xi_2^\star+ r \rho \cos(\psi-\varphi), 
\end{align*}
and that, therefore, the leading asymptotic contribution of $\bdxi_{\text{SD}1}^{(F)}$ is given by 
\begin{align*}
  u_{\text{SD} 1}(x_1,x_2,0^+)  & =  \mathcal{A}e^{- i(x_1 \xi_1^\star+x_2 \xi_2^\star)}
  \iint_{\Psi (\tmmathbf{\Gamma})} \frac{(\cos (\psi))^{1 / 2} e^{- i \rho r
  \cos (\psi - \varphi)}}{(\sin(\vartheta_1-\psi))^{1 / 2}
  \sin (\psi)} \mathd \rho \wedge \mathd \psi,
\end{align*}
where we introduced the notation 
\[
\mathcal{A} = \frac{\sqrt{k'_1 + k_1}}{4\pi^2\sqrt{2k}(-k_1'+k_1)}.
\] 

Our aim is hence to evaluate this double integral. It can be shown that the surface of integration $\Psi (\tmmathbf{\Gamma})$ can
be can be continuously deformed (without hitting any singularities) to the
surface $\left\{ (\rho, \psi) \in \mathbb{C}^2, \psi \in \Upsilon \infixand
\rho \in \lambda_{\psi} \right\}$. The contour $\Upsilon$, illustrated in \figurename~\ref{fig:Phi_contour} (left), is close to the real segment $[0, 2 \pi]$. The way that $\Upsilon$ bypasses the singular points in the $\psi$ plane provides the required compatibility with the bridge configuration around the triple crossing.

The $\lambda_{\psi}$ form a continuous family of contours in
$\mathbb{C}$ that depend on $\text{Re}[\psi]$ and vary between $0$ and $\infty$, as shown in \figurename~\ref{fig:Phi_contour} (right). It is chosen to ensure exponential attenuation of the integrand as $\rho$ tends to $\infty$.  Note that when $\text{Re}[\psi]$ is equal to $\varphi+\pi/2$ or $\varphi+3\pi/2$, the contour $\lambda_\psi$ has to be purely real. However, in that case, the exponential decay is provided by the imaginary part of $\psi$. This is why, even though the points $\varphi+\pi/2$ and $\varphi+3\pi/2$ are not singular, they are chosen to be bypassed by $\Upsilon$ as in \figurename~\ref{fig:Phi_contour} (left).

\begin{figure}[h]
  \centering{\includegraphics[width=0.54\textwidth]{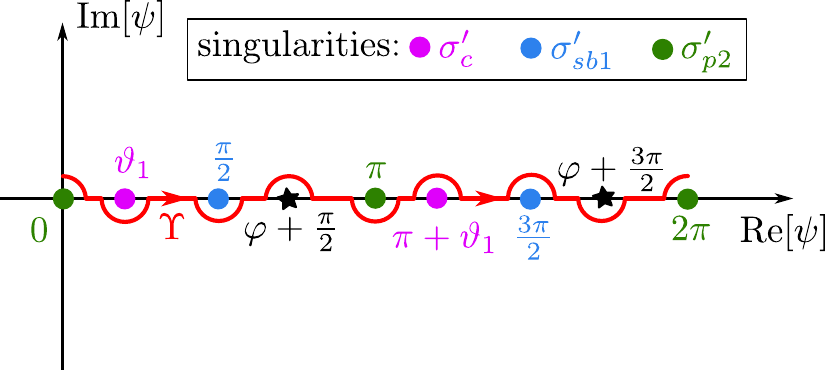} \includegraphics[width=0.44\textwidth]{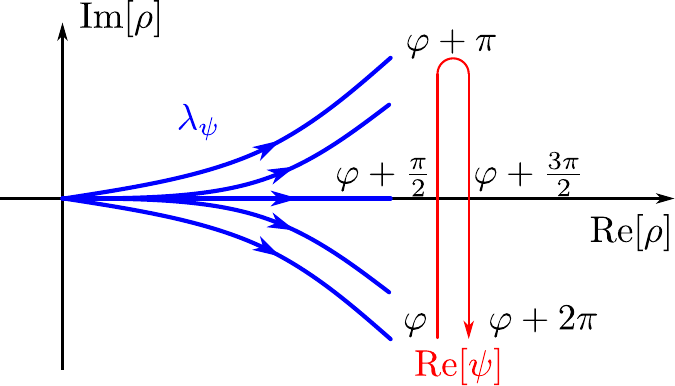}}
  \caption{The contour $\Upsilon$ (left) and the family of contours $\lambda_\psi$ (right)}
\label{fig:Phi_contour}
\end{figure}

%Note that for two values of angular variable $\psi = \varphi + \pi/2$, $\psi = \varphi + 3\pi/2$ the radial variable $\rho$ is real, and attenuation is provided by the imaginary part of $\psi$ (see bypasses around these points in  \figurename~\ref{fig:Phi_contour}).

%\begin{figure}[h]
%  \centering{\includegraphics[width=0.5\textwidth]{pictures/Contour_lambda_Psi.pdf}}
%  \caption{Family integration contours $\lambda_\psi$}
%\label{fig:lambda_contours}
%\end{figure}

   The integral can hence be rewritten in a sequential form as
follows:
\begin{eqnarray*}
  u_{\text{SD} 1}(x_1,x_2,+0)  & \approx & \mathcal{A}e^{- i(x_1 \xi_1^\star+x_2 \xi_2^\star)}
  \int_{\psi \in \Upsilon} \frac{(\cos (\psi))^{1 / 2} }{(\sin(\vartheta_1-\psi))^{1 / 2} \sin (\psi)}  \left( \int_{\rho \in
  \lambda_{\psi}} e^{- i \rho r \cos (\psi - \varphi)} \mathd \rho \right)
  \mathd \psi .
\end{eqnarray*}
Given the properties of $\lambda_{\psi}$, the inner integral can be evaluated directly to
\begin{eqnarray*}
  \int_{\rho \in \lambda_{\psi}} e^{- i \rho r \cos (\psi - \varphi)} \mathd
  \rho & = & \frac{- i}{r \cos (\psi - \varphi)},
\end{eqnarray*}
and $u_{\text{SD} 1}$ is now given by a single integral
\begin{align}
  u_{\text{SD} 1}(x_1,x_2,0^+)  & \approx & \frac{- i\mathcal{A}e^{- i(x_1 \xi_1^\star+x_2 \xi_2^\star)}}{r} \int_{\Upsilon} \frac{(\cos (\psi))^{1 / 2} }{(\sin(\vartheta_1-\psi))^{1 / 2} \sin (\psi) \cos (\psi -
  \varphi)} \mathd \psi .
\label{eq:integralItriplecrossintermediate}
\end{align}
%At this stage, since $(k_1')^2 + (k_2)^2 = k^2$, we can introduce an angle
%$\vartheta_1$ such that
%\begin{eqnarray*}
%  k_1' = k \sin (\vartheta_1) & \tmop{and} & k_2 = - k \cos (\vartheta_1)
%\end{eqnarray*}
%and the integral becomes
%\begin{eqnarray}
%  u_{\text{SD} 1}(x_1,x_2,+0)  & \approx & \frac{- i\mathcal{A}e^{- i\tmmathbf{x} \cdot
%  \tmmathbf{\xi}^{\star}}}{r \sqrt{2 k}} \int_{\Upsilon} \frac{(\cos
%  (\psi))^{1 / 2} }{(\sin (\theta_1) \cos (\psi) - \cos (\vartheta_1) \sin
%  (\psi))^{1 / 2} \sin (\psi) \cos (\psi - \varphi)} \mathd \psi \nonumber\\
%  & \approx & \frac{- i\mathcal{A}e^{- i\tmmathbf{x} \cdot
%  \tmmathbf{\xi}^{\star}}}{r \sqrt{2 k}} \int_{\Upsilon} \frac{(\cos
%  (\psi))^{1 / 2} }{(\sin (\vartheta_1 - \psi))^{1 / 2} \sin (\psi) \cos (\psi -
%  \varphi)} \mathd \psi.  \label{eq:integralItriplecrossintermediate}
%\end{eqnarray}
Fortunately, this integral can be evaluated exactly (see~Appendix~\ref{app:triple_crossing} for details): 
\begin{equation*}
  \int_{\Upsilon} \frac{(\cos (\psi))^{1 / 2} }{(\sin (\vartheta_1 - \psi))^{1 /
  2} \sin (\psi) \cos (\psi - \varphi)} \mathd \psi  =  \frac{4\pi \sqrt{-\sin(\varphi)}}{\cos(\varphi)\sqrt{\cos(\vartheta_1-\varphi)}}\mathcal{H}(-\sin(\varphi))\mathcal{H}(\cos(\vartheta_1-\varphi)).
\end{equation*}
It leads to the following expression for the secondary diffracted field at $x_3 = 0$:
\begin{equation}
  u_{\text{SD} 1} (x_1,x_2,0^+) =  \frac{-i\sqrt{k'_1+k_1}\sqrt{-x_2}e^{i(k'_1x_1 + k_2x_2)}}{\sqrt{2}\pi(k_1-k_1')x_1\sqrt{k'_1x_2-k_2x_1}}\mathcal{H} (-x_2) \mathcal{H}(k'_1x_2-k_2x_1). 
  \label{eq:usd1formula}
\end{equation}
The asymptotic contribution resulting from the other triple crossing $\tmmathbf{\xi}_{\text{SD} 2}^{(F)}$
can be obtained similarly, or just by swapping the subscripts 1 and 2, and reads
\begin{equation}
  u_{\text{SD} 2} (x_1,x_2,0^+) =  \frac{-i\sqrt{k'_2+k_2}\sqrt{-x_1}e^{i(k_1x_1+k'_2x_2 )}}{\sqrt{2}\pi(k_2-k_2')x_2\sqrt{k'_2x_1-k_1x_2}}\mathcal{H} (-x_1) \mathcal{H}(k'_2x_1-k_1x_2). 
  \label{eq:usd2formula}
\end{equation}

\begin{rema}
  While doing this work, with the aim of double checking our result, we found a way of recovering this formula by using the formalism of \cite{Shanin2012}.
\end{rema}

\subsection{General asymptotic expansion of $u$}

We have now considered all the contributing points and we can therefore
reconstruct the far-field asymptotic approximation of $u$ as follows:
\begin{eqnarray*}
  u (\tmmathbf{x}) & \underset{| \tmmathbf{x} | \rightarrow \infty}{\approx} &
  u_{\tmop{RW}} (\tmmathbf{x}) + u_{\text{PD} 1}  (\tmmathbf{x}) +
  u_{\text{PD} 2} (\tmmathbf{x}) + u_{\tmop{SW}} (\tmmathbf{x}) + u_{\text{SD}
  1} (\tmmathbf{x}) + u_{\text{SD} 2} (\tmmathbf{x}).
\end{eqnarray*}
It agrees, in principle, with former results given for the field in {\cite{Assier2012,Shanin2012,Lyalinov2013}}. This expansion is valid for $x_3>0$ and for $x_3=0$. When $x_3>0$, the components are given by equations (\ref{eq:uPWsimple}), (\ref{eq:uPD1simple}), (\ref{eq:uPD2simple}), (\ref{eq:uSWsimple}), (\ref{eq:usd1formula_x3neq0}),(\ref{eq:usd2formula_x3neq0}). The first four components are continuous as $x_3\to0$, however, more work is required for the last two secondary diffracted waves, for which the formulae  (\ref{eq:usd1formula}) and (\ref{eq:usd2formula}) should be used when $x_3=0$. A schematic illustration of these wave components is given in \figurename~\ref{fig:far-field-complicated}.

\begin{figure}[h]
	\centering{
		\includegraphics[height=0.4\textwidth]{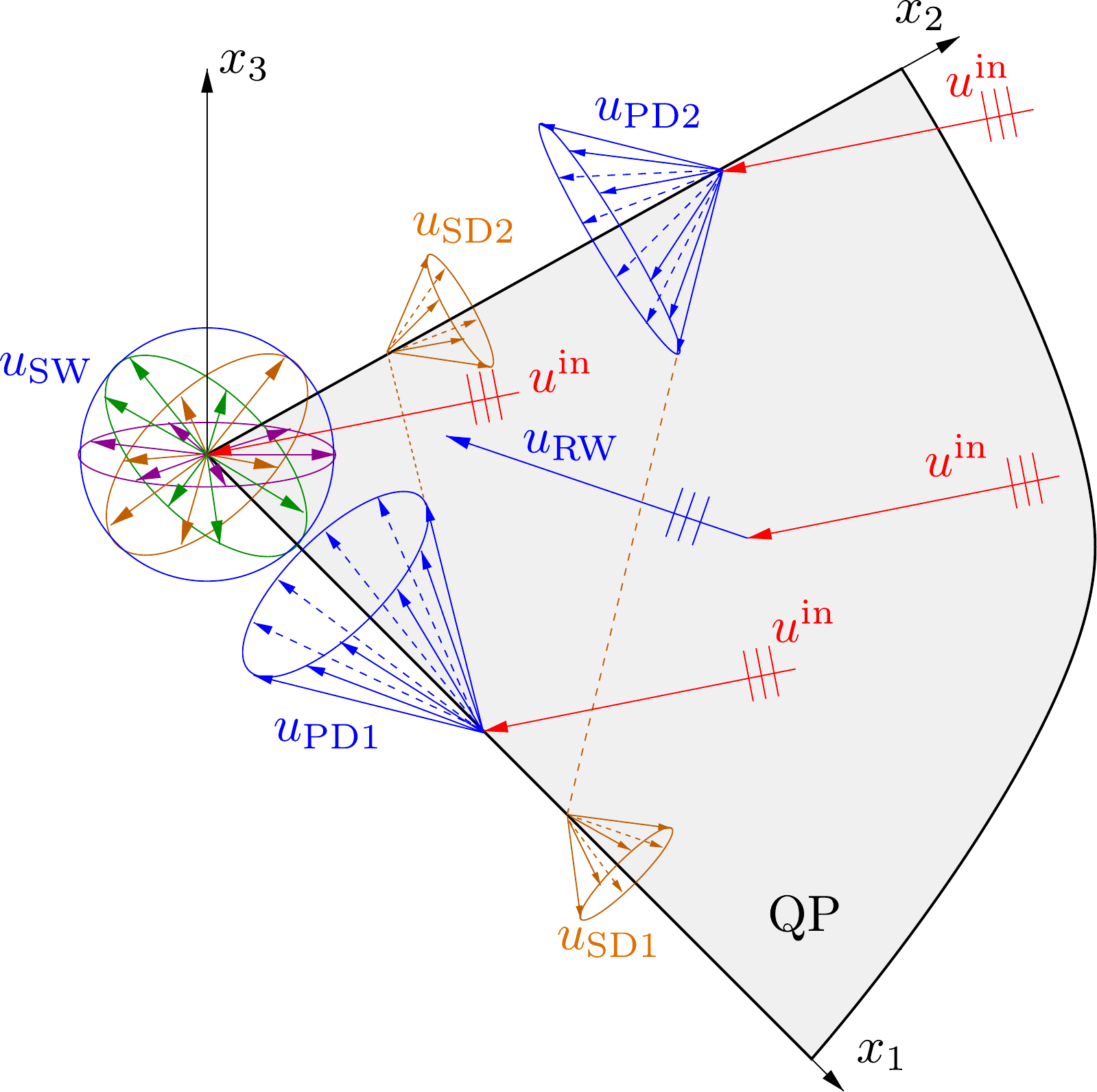}
	}
	\caption{Illustration of the far-field asymptotic components in complicated case (adapted from~\cite{Assier2012}).}
	\label{fig:far-field-complicated}
\end{figure}

%These expansions agree in principle with former results given for the field in
%{\cite{Assier2012,Shanin2012,Lyalinov2013}}, and additionally we have obtained the formulae  (\ref{eq:usd1formula}) and (\ref{eq:usd2formula}) for
%$u_{\text{SD} 1} (x_1,x_2,0^+)$ and $u_{\text{SD} 2} (x_1,x_2,0^+)$.

\section{Conclusion}

We have shown that it was possible to apply the methodology developed in
{\cite{Part6A}} to the challenging problem of wave diffraction by a quarter-plane.
We used this methodology, together with the analytical continuation formulae
derived in {\cite{Assier2018a}} to recover the far-field asymptotics of the
problem at hand with only a modest amount of algebra. We recovered all the
known results on the far-field asymptotics of the wave field in two cases, the simple and the complicate case. All
the far-field wave components were obtained by a direct application of the
results of {\cite{Part6A}}. The only formulae that could not be obtained in such a way concern the approximation of the secondary diffracted waves (only occurring in the complicated case) on the plane $x_3=0$. Those were found to correspond to a triple crossing of
singularities. We however managed to adapt our argument and recovered a
far-field component, leading to new formulae that had not been obtained
before. Though we have chosen, for brevity, not to deal in details with the normal derivative of the field, its asymptotics can be obtained very similarly by studying (\ref{e:sh010}) instead of (\ref{e:sh009})$|_{x_3 = 0}$. In that case, as illustrated in \figurename~\ref{fig:complicated-case} (right), we would not have to deal with triple crossings, so the components can be obtained directly by applying the methodology of \cite{Part6A}.

%%% Bibliography %%%
%\small
\bibliography{biblio}
\bibliographystyle{unsrt}

\appendix
\section{Proof of additive
crossing}\label{app:proof-of-additive-crossing}

Though we have concluded in \cite{Assier2018a} that the crossing $(-k,-k)$ must be additive, we provide here an additional proof of this fact based solely on the analytic continuation formulae. Consider the formula (\ref{eq:third-anal-formula-W}) and take $\xi_1\in H^-$ and $\xi_2 \in h^-$, with $\xi_2$ being on the left of $h^-$, as illustrated in \figurename~\ref{fig01a}. Since $\xi_2$ is on $P$, the integral along $P$ is understood in the sense of principle value. 
%According to the concept of the additive crossing, 
Let us also introduce the path $\delta_1'$ (resp. $\delta_2'$),  which is a loop starting and ending at $\xi_1$ (resp. $\xi_2$) bypassing $-k$ anticlockwise (see \figurename~\ref{fig01a}).  

%%%%%%%%%%%%%%%%%%%%%%%
\begin{figure}[h]
\centering
\includegraphics[width = 0.9\textwidth]{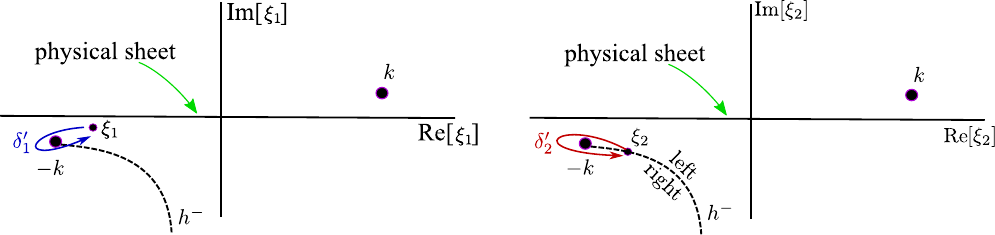}
\caption{The bypasses $\delta_1'$ and $\delta_2'$}
\label{fig01a}
\end{figure}
%%%%%%%%%%%%%%%%%%%%%%%
Rewrite (\ref{eq:third-anal-formula-W}) as
\begin{equation}
W (\xi_1 , \xi_2) = \frac{\tilde J_1(\xi_1 , \xi_2)}{4 \pi^2}
+ 
\frac{i \gamma(\xi_1 , \xi_2) \gamma(\xi_1 , k_2) \gamma(k_2 , k_1)}{(\xi_1 + k_1) (\xi_2 + k_2) \gamma(k_2 , - \xi_1)}
, 
\label{e6015b}
\end{equation}
\begin{equation}
	\tilde J_1(\xi_1 , \xi_2) \equiv
	\gamma(\xi_1 , \xi_2) 
	\,   \int \limits_{-\infty - i\epsilon }^{\infty - i\epsilon} \left( 	\int \limits_{P}  \frac{\gamma(\xi_1 , - \xi_2') K (\xi_1' , \xi_2') W(\xi_1' , \xi_2')
	}{ (\xi_2' - \xi_2)} d\xi_2' \right) \frac{d\xi_1'}{(\xi_1' - \xi_1)}
	\label{e6015c}
\end{equation}
where $\epsilon$ is a positive value such that $0 < \epsilon < -{\rm Im}[\xi_1]$.
Note that we can make this deformation due to the position of~$P$.
We introduce four values of $W (\xi_1 , \xi_2)$, they are $W (\xi_1 , \xi_2)$, $W_{\delta_1'} (\xi_1, \xi_2)$,
$W_{\delta_2'} (\xi_1, \xi_2)$, $W_{\delta_1' \delta_2'} (\xi_1, \xi_2)$, 
where $W (\xi_1 , \xi_2)$ is the ``initial'' value,
$W_{\delta_1'} (\xi_1 , \xi_2)$ is obtained from  $W (\xi_1 , \xi_2)$  by continuation along $\delta_1'$ in the $\xi_1$-plane, 
$W_{\delta_2'} (\xi_1 , \xi_2)$ is obtained from  $W (\xi_1 , \xi_2)$  by continuation along $\delta_2'$ in the $\xi_2$-plane, 
$W_{\delta_1'\delta_2'} (\xi_1  \xi_2)$ is obtained from  $W_{\delta_1'} (\xi_1 , \xi_2)$  
by continuation along $\delta_2'$ in the $\xi_2$-plane. To prove that the crossing $\bdxi^\star=(-k,-k)$ is additive for $W$, it is enough to show that 
\begin{equation}
W (\xi_1 , \xi_2) +W_{\delta_1' \delta_2'} (\xi_1 , \xi_2) = 
W_{\delta_1'} (\xi_1 , \xi_2) + W_{\delta_2'} (\xi_1 , \xi_2),
\label{e6015a}
\end{equation}
as $\xi_1$ approaches $h^-$.
Indeed, if (\ref{e6015a}) is satisfied then $W$ can be presented in the form (\ref{e:sh026}). This is due to {\color{magenta}Proposition 4.12} of \cite{Part6A} (which was first proved in Section 4.2 of \cite{Assier2018a}). 
Because it only contains usual functions, one can check directly that the second term of (\ref{e6015b}) possesses the additive crossing property for $\bdxi^\star=(-k,-k)$ by showing that it satisfies a relationship akin to (\ref{e6015a}). Thus, we only have to prove it for the term $\tilde J_1$.
The integral (\ref{e6015c}) cannot be understood literally, since the point $\xi_2$ sits exactly on $P$. Rewriting the inner integral in (\ref{e6015c}) as an integral along the left part of $h^-$ only, and using the Sokhotski-Plemelj formula, we can rewrite $\tilde J_1$ as:
%Let us rewrite the
%integral for $\tilde J_1$, $\tilde J_{1,\delta_1'}$, $\tilde J_{1,\delta_2'}$, $\tilde J_{1,\delta_1' \delta_2'}$ (the indices are assigned as above) using the Sokhotski-Plemelj formula:
\begin{align}
\tilde J_1(\xi_1 , \xi_2) &\equiv
\gamma(\xi_1 , \xi_2) 
\int \limits_{-\infty - i\epsilon }^{\infty - i\epsilon} d\xi_1' \, {\rm V.P.} \!\!\! \int \limits_{P} d\xi_2' 
\frac{\gamma(\xi_1 , - \xi_2') K (\xi_1' , \xi_2') W(\xi_1' , \xi_2')
}{ (\xi_1' - \xi_1)(\xi_2' - \xi_2)} \nonumber \\ 
&+ 
\frac{\pi i }{ K (\xi_1 , \xi_2)}
\int \limits_{-\infty - i\epsilon }^{\infty - i\epsilon} d\xi_1'
\frac{K (\xi_1' , \xi_2) 
[W(\xi_1' , \xi_2) - W_{\delta_2'}(\xi_1' , \xi_2)]
}{
\xi_1' - \xi_1
},
\label{e6015d}
\end{align}
where V.P.\ denotes the principal value of a singular integral, defined in a usual way, 
i.e.\ we consider separately small arcs bypassing the point $\xi_2' = \xi_2$ (see \figurename~\ref{fig01b}).
From this we can then obtain the following expressions for $\tilde J_{1,\delta_1'}$, $\tilde J_{1,\delta_2'}$, $\tilde J_{1,\delta_1' \delta_2'}$ (the indices are assigned as they were above for $W$):

%%%%%%%%%%%%%%%%%%%%%%%
\begin{figure}[h]
\centering
\includegraphics[width = 0.6\textwidth]{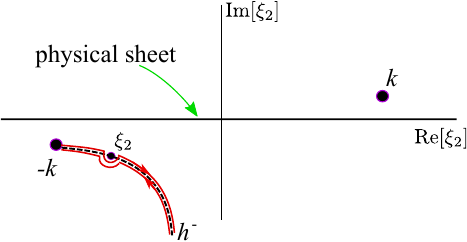}
\caption{The contour $P$ and the small arc deformations used for the V.P. integrals}
\label{fig01b}
\end{figure}
%%%%%%%%%%%%%%%%%%%%%%%
\begin{align}
\tilde J_{1,\delta_1'} (\xi_1 , \xi_2) &\equiv
\gamma' (\xi_1 , \xi_2)
\int \limits_{-\infty - i\epsilon }^{\infty - i\epsilon} d\xi_1' \, {\rm V.P.} \!\!\! \int \limits_{P} d\xi_2' 
\frac{ \gamma'(\xi_1 , - \xi_2') K (\xi_1' , \xi_2') W(\xi_1' , \xi_2')
}{ (\xi_1' - \xi_1)(\xi_2' - \xi_2)} \nonumber \\
&+ 
\frac{\pi i }{ K (\xi_1 , \xi_2)}
\int \limits_{-\infty - i\epsilon }^{\infty - i\epsilon} d\xi_1'
\frac{K (\xi_1' , \xi_2) 
[W(\xi_1' , \xi_2) - W_{\delta_2'}(\xi_1' , \xi_2)]
}{
\xi_1' - \xi_1
} , 
\label{e6015f}
\end{align}
\begin{align}
\tilde J_{1,\delta_2'} (\xi_1 , \xi_2) &\equiv
\gamma(\xi_1 , \xi_2) 
\int \limits_{-\infty - i\epsilon }^{\infty - i\epsilon} d\xi_1' \, {\rm V.P.}\!\!\! \int \limits_{P} d\xi_2' 
\frac{\gamma(\xi_1 , - \xi_2') K (\xi_1' , \xi_2') W(\xi_1' , \xi_2')
}{ (\xi_1' - \xi_1)(\xi_2' - \xi_2)} \nonumber \\
&- 
\frac{\pi i }{ K (\xi_1 , \xi_2)}
\int \limits_{-\infty - i\epsilon }^{\infty - i\epsilon} d\xi_1'
\frac{K (\xi_1' , \xi_2) 
[W(\xi_1' , \xi_2) - W_{\delta_2'}(\xi_1' , \xi_2)]
}{
\xi_1' - \xi_1
} , 
\label{e6015e}
\end{align}
\begin{align}
\tilde J_{1,\delta_1' \delta_2'} (\xi_1 , \xi_2) &\equiv
\gamma'(\xi_1, \xi_2)
\int \limits_{-\infty - i\epsilon }^{\infty - i\epsilon} d\xi_1' \, {\rm V.P.} \!\!\! \int \limits_{P} d\xi_2' 
\frac{ \gamma'(\xi_1 , -\xi_2') K (\xi_1' , \xi_2') W(\xi_1' , \xi_2')
}{ (\xi_1' - \xi_1)(\xi_2' - \xi_2)} \nonumber \\
&- 
\frac{\pi i }{ K (\xi_1 , \xi_2)}
\int \limits_{-\infty - i\epsilon }^{\infty - i\epsilon} d\xi_1'
\frac{K (\xi_1' , \xi_2) 
[W(\xi_1' , \xi_2) - W_{\delta_2'}(\xi_1' , \xi_2)]
}{
\xi_1' - \xi_1
} , 
\label{e6015g}
\end{align}
where we have defined
\[
\gamma'(\xi_1 , \xi_2) \equiv \sqrt{- \sqrt{k^2 - \xi_1^2} + \xi_2}. 
\]
The additive crossing relation 
\begin{equation}
\tilde J_1 (\xi_1 , \xi_2) + \tilde J_{1,\delta_1' \delta_2'} (\xi_1 , \xi_2)
= 
\tilde J_{1,\delta_1'} (\xi_1 , \xi_2) + \tilde J_{1,\delta_2'} (\xi_1 , \xi_2)
\label{e6015h}
\end{equation}
can now be checked directly.
Note that the principal value part of $\tilde J_1$ does not change after the bypass $\delta_2'$, while
the residue part does not change after the bypass~$\delta_1'$. We have therefore proven that (\ref{e6015a}) is true on $H^- \times h^-$. Then, by analytical continuation in the $\xi_1$ plane, it is true on $h^- \times h^-$, as required, and the crossing is additive.

\section{Evaluation of  (\ref{eq:integralItriplecrossintermediate})}
\label{app:triple_crossing} 
Let us evaluate the integral:
\begin{equation}
  I  = \int_{\Upsilon} \text{Int}(\psi) \mathd \psi, \text{ where } \text{Int}(\psi)=\frac{(\cos
  	(\psi))^{1 / 2} }{(\sin (\vartheta_1 - \psi))^{1 / 2} \sin (\psi) \cos (\psi -
  	\varphi)} \cdot
\end{equation}
To do this, we will consider three cases, depending on the location of $\varphi+\pi/2$:
\begin{align}
	\text{Case 1: } 0<\varphi+\pi/2<\theta_1,\quad \text{Case 2: } \theta_1<\varphi+\pi/2<\pi/2, \quad \text{Case 3: } \pi/2<\varphi+\pi/2<\pi.
\end{align}
The contour $\Upsilon$ is shown for each case in \figurename~\ref{fig:Phi_contour_def} (top row). On this figure, the dotted lines represent branch cuts of the integrand $\text{Int}(\psi)$. In each case, we start by deforming the contour $\Upsilon$ as described in \figurename~\ref{fig:Phi_contour_def} (middle row). We note that in the Cases 2 \& 3, two poles are picked up in the process. It can be shown directly that the integrand  $\text{Int}(\psi)\to 0$ exponentially as $\text{Im}[\psi]\to\pm\infty$. Hence, we can push the deformation process further and discard the horizontal parts of the contours. The resulting contours for each case are displayed in \figurename~\ref{fig:Phi_contour_def} (bottom row). They consist of two infinite vertical lines with opposite orientation, and, in Cases 2 \& 3, some encircled poles. The two lines are at a distance $\pi$ from each other and therefore $\text{Int}(\psi)$ takes the same values on each line. The part of the integral corresponding to these lines can hence be discarded.
% First, one can check directly that the antiderivative of the integrand is given by:
%\begin{align*}
%%  &\int \frac{(\cos (\psi))^{1 / 2} }{(\sin (\vartheta_1 - \psi))^{1 / 2} \sin
%%  (\psi) \cos (\psi - \varphi)} \mathd \psi  \\ = 
%2 \sec (\varphi) \sqrt{\cos
%  (\psi)} &\left( \frac{\sqrt{\sin (\varphi)}  \sqrt{\sec (\psi)} \tanh^{- 1}
%  \left( \frac{\sqrt{\sin (\varphi)}  \sqrt{\sec (\psi)}  \sqrt{\sin (\vartheta_1
%  - \psi)}}{\sqrt{\cos (\vartheta_1 - \varphi)}} \right)}{\sqrt{\cos (\vartheta_1 -
%  \varphi)}} \right.\\
%  & \quad -  \left. \frac{\sqrt{\sec (\psi) \sin (\vartheta_1 - \psi)} \tanh^{- 1}
%  \left( \frac{\sqrt{\sin (\vartheta_1) - \cos (\theta_1) \tan
%  (\psi)}}{\sqrt{\sin (\vartheta_1)}} \right)}{\sqrt{\sin (\vartheta_1)}  \sqrt{\sin
%  (\vartheta_1 - \psi)}} \right)+\text{cst}. 
%\end{align*}
%Second, let us represent $I$ as a sum of two integrals:
%\begin{equation}
%I = \frac{1}{2}(I_{\Upsilon_1} + I_{\Upsilon_2}),
%\end{equation}
%where $I_{\Upsilon_1}$ and  $I_{\Upsilon_2}$ have the same integrand as $I$ but with contours of integration deformed into $\Upsilon_1$ and $\Upsilon_2$ correspondingly. The deformed contours are presented in \figurename~\ref{fig:Phi_contour_def}. $\Upsilon_1$ (resp.\ $\Upsilon_1$) is decomposed into a straight line $\gamma_u$ (resp.\ $\gamma_d$), a contour around a branch cut $\gamma_\ell$ (resp. $\gamma_r$) and two pole-encircling contours.
\begin{figure}[h!]
  \centering{\includegraphics[width=1\textwidth]{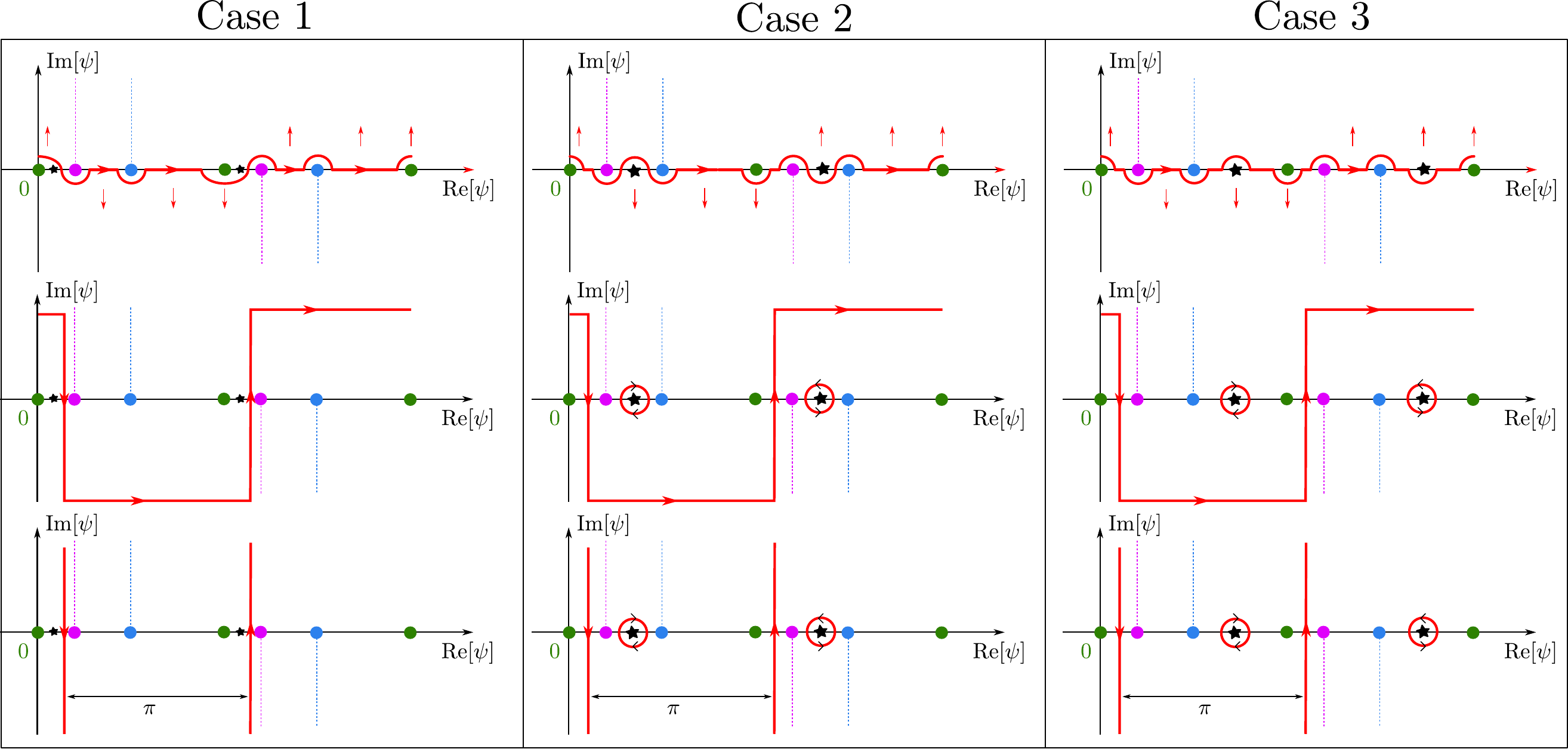}}
  \caption{Process of contour deformation of $\Upsilon$ in the three cases}
\label{fig:Phi_contour_def}
\end{figure}

We are left with the following result:
\begin{align}
	I & =  \left\{ \begin{array}{cl}
		0 & \text{in Case 1}\\
		- 2 i \pi \tmop{Res} \left[ \tmop{Int}(\psi), \psi = \varphi + \frac{\pi}{2}
		\right] + 2 i \pi \tmop{Res} \left[ \tmop{Int}(\psi), \psi = \varphi + \frac{3
			\pi}{2} \right] & \text{in Cases 2 \& 3}
	\end{array} \right.
\label{eq:raphstep2appB}
\end{align}
Using the fact that in Case 2 we have $\sin(\varphi)<0$ and $\cos(\theta_1-\varphi)>0$, while in Case 3 we have $\sin(\varphi)>0$ and $\cos(\theta_1-\varphi)>0$, we find that

\begin{alignat*}{2}
	\tmop{Res} \left[ \tmop{Int}(\psi), \psi = \varphi + \frac{\pi}{2} \right] &=
	\frac{- (- \sin (\varphi))^{1 / 2}}{(- \cos (\vartheta_1 - \varphi))^{1 / 2}
		\cos (\varphi) }  &&=  \left\{ \begin{array}{cc}
		\frac{i \sqrt{- \sin (\varphi)}}{\cos
			(\varphi)\sqrt{\cos (\vartheta_1 - \varphi)}  } &\text{in Case 2}\\
		\frac{- \sqrt{\sin (\varphi)}}{ \cos
			(\varphi)\sqrt{\cos (\vartheta_1 - \varphi)} } &\text{in Case 3}
	\end{array} \right. ,\\
	\tmop{Res} \left[ \tmop{Int}(\psi), \psi = \varphi + \frac{3 \pi}{2} \right] &=
	\frac{- (\sin (\varphi))^{1 / 2}}{(\cos (\vartheta_1 - \varphi))^{1 / 2}
		\cos (\varphi)}  &&=  \left\{ \begin{array}{cc}
		\frac{- i \sqrt{- \sin (\varphi)}}{\cos (\varphi)\sqrt{\cos (\vartheta_1 - \varphi)}
			} & \text{in Case 2}\\
		\frac{- \sqrt{\sin (\varphi)}}{\cos
			(\varphi)\sqrt{\cos (\vartheta_1 - \varphi)}  } & \text{in Case 3}
	\end{array} \right. .
\end{alignat*}
Inputting this into (\ref{eq:raphstep2appB}), we find that 
\begin{align*}
	I & =  \left\{ \begin{array}{cl}
		0 & \text{in Cases 1 \& 3}\\
		\frac{4 \pi \sqrt{- \sin (\varphi)}}{\sqrt{\cos (\vartheta_1 - \varphi)}
			\cos (\varphi) } & \text{in Case 2}
	\end{array} \right. ,
\end{align*}
which can be summarised by
\begin{equation}
	I  =  \frac{4\pi \sqrt{-\sin(\varphi)}}{\cos(\varphi)\sqrt{\cos(\vartheta_1-\varphi)}}\mathcal{H}(-\sin(\varphi))\mathcal{H}(\cos(\vartheta_1-\varphi)),
\end{equation}
as required, since, in Case 1, we have $\cos(\vartheta_1-\varphi)<0$.

\end{document}